\documentclass[12pt]{article}  
\input amssym.def  
\input amssym.tex

\usepackage[frame,cmtip,curve,arrow,matrix,line,graph]{xy}

\usepackage{euscript}  
  
\begin{document}

\newcommand{\nc}{\newcommand}  
  
\newcommand{\colvec}[2]{\left  ( \begin{array}{cc} #1  \\  
     #2  \end{array} \right ) }  
  
\newcommand{\Tr}{\,{\rm Tr}\,}  
\newcommand{\End}{\,{\rm End}\,}  
\newcommand{\Hom}{\,{\rm Hom}\,}  
  
\newcommand{\Ker}{ \,{\rm Ker} \,}  
  
\newcommand{\bla}{\phantom{bbbbb}}  
\newcommand{\onebl}{\phantom{a} }  
\newcommand{\eqdef}{\;\: {\stackrel{ {\rm def} }{=} } \;\:}  
\newcommand{\sign}{\: {\rm sign}\: }  
\newcommand{\sgn}{ \:{\rm sgn}\:}  
\newcommand{\half}{ {\frac{1}{2} } }  
\newcommand{\vol}{ \,{\rm vol}\, }

%
  
\newcommand{\beq}{\begin{equation}}  
\newcommand{\eeq}{\end{equation}}  
\newcommand{\beqst}{\begin{equation*}}  
\newcommand{\eeqst}{\end{equation*}}  
\newcommand{\barr}{\begin{array}}  
\newcommand{\earr}{\end{array}}  
\newcommand{\beqar}{\begin{eqnarray}}  
\newcommand{\eeqar}{\end{eqnarray}}  
\newtheorem{theorem}{Theorem}
\newtheorem{corollary}[theorem]{Corollary}  
\newtheorem{lemma}[theorem]{Lemma}  
\newtheorem{prop}[theorem]{Proposition}  
\newtheorem{proposition}[theorem]{Proposition}  
\newtheorem{definition}[theorem]{Definition}  
\newtheorem{examplit}[theorem]{Example}  
\newtheorem{remit}[theorem]{Remark}  
\newtheorem{conjecture}[theorem]{Conjecture}  
\newcommand{\matr}[4]{\left \lbrack \begin{array}{cc} #1 & #2 \\  
     #3 & #4 \end{array} \right \rbrack}

\newenvironment{example}{\begin{examplit}\rm}{\end{examplit}}  
\newenvironment{rem}{\begin{remit}\rm}{\end{remit}}

\newcommand{\aff}{{\Bbb A }}  
\newcommand{\RR}{{\Bbb R }}  
\newcommand{\CC}{{\Bbb C }}  
\newcommand{\ZZ}{{\Bbb Z }}  
\newcommand{\PP}{ {\Bbb P } }  
\newcommand{\QQ}{{\Bbb Q }}  
\newcommand{\UU}{{\Bbb U }}

\newcommand{\cala}{{\mbox{$\mathcal A$}}}  
\newcommand{\calb}{{\mbox{$\mathcal B$}}}  
\newcommand{\calc}{{\mbox{$\mathcal C$}}}  
\newcommand{\cald}{{\mbox{$\mathcal D$}}}  
\newcommand{\cale}{{\mbox{$\mathcal E$}}}  
\newcommand{\calf}{{\mbox{$\mathcal F$}}}  
\newcommand{\calg}{{\mbox{$\mathcal G$}}}  
\newcommand{\calh}{{\mbox{$\mathcal H$}}}  
\newcommand{\cali}{{\mbox{$\mathcal I$}}}  
\newcommand{\calj}{{\mbox{$\mathcal J$}}}  
\newcommand{\calk}{{\mbox{$\mathcal K$}}}  
\newcommand{\call}{{\mbox{$\mathcal L$}}}  
\newcommand{\calm}{{\mbox{$\mathcal M$}}}  
\newcommand{\caln}{{\mbox{$\mathcal N$}}}  
\newcommand{\calo}{{\mbox{$\mathcal O$}}}  
\newcommand{\calp}{{\mbox{$\mathcal P$}}}  
\newcommand{\calq}{{\mbox{$\mathcal Q$}}}  
\newcommand{\calr}{{\mbox{$\mathcal R$}}}  
\newcommand{\cals}{{\mbox{$\mathcal S$}}}  
\newcommand{\calt}{{\mbox{$\mathcal T$}}}  
\newcommand{\calu}{{\mbox{$\mathcal U$}}}  
\newcommand{\calv}{{\mbox{$\mathcal V$}}}  
\newcommand{\calw}{{\mbox{$\mathcal W$}}}  
\newcommand{\calx}{{\mbox{$\mathcal X$}}}  
\newcommand{\caly}{{\mbox{$\mathcal Y$}}}  
\newcommand{\calz}{{\mbox{$\mathcal Z$}}}

\def\a{\alpha}  
\def\b{\beta}  
\def\g{\gamma}  
\def\d{\delta}  
\def\e{\epsilon}  
\def\z{\zeta}  
\def\h{\eta}  
\def\t{\theta}  
\def\k{\kappa}  
\def\l{\lambda}  
\def\m{\mu}  
\def\n{\nu}  
\def\x{\xi}  
\def\p{\pi}  
\def\r{\rho}  
\def\s{\sigma}  
\def\ta{\tau}  
\def\u{\upsilon}  
\def\ph{\phi}  
\def\c{\chi}  
\def\ps{\psi}  
\def\o{\omega}  
  
\def\G{\Gamma}  
\def\D{\Delta}  
\def\T{\Theta}  
\def\L{\Lambda}  
\def\X{\Xi}  
\def\P{\Pi}  
\def\S{\Sigma}  
\def\U{\Upsilon}  
\def\Ph{\Phi}  
\def\Ps{\Psi}  
\def\O{\Omega}

\nc{\Proof}{\noindent{\em Proof:} }   
\nc{\Lie}{ {\rm Lie} }   
\nc{\liek}{{ \bf k} }  
\nc{\lieks}{\liek^*}  
\nc{\hk}{H^*_K}  
\nc{\hT}{H^*_T}  
\nc{\xvec}{X}  
\nc{\hfeta}{h_F^{\eta}}

\nc{\lieksa}{\lieks_{{\rm ad}} }  
\nc{\lieka}{\liek_{{\rm ad}} }  
\nc{\liekas}{(\liek_s)_{{\rm ad}} }  
\nc{\lieksas}{(\lieks_s)_{{\rm ad}} }  
\nc{\liekaq}{(\liek/\liek_s)_{{\rm ad}} }  
\nc{\lieksaq}{(\lieks/\lieks_s)_{{\rm ad}} }  
\nc{\conv}{{\rm Conv}}  
  
\nc{\proj}{\pi}  
\nc{\zloc}{\mu^{-1}(0)}   
\nc{\nf}{\nu_F}   
\nc{\lx}{\l(X) }   
\nc{\rlx}{\frac{d \l(X)}{\l(X) }  }   
\nc{\srlx}{ {d \l(X)} / {\l(X) }  }   
\nc{\intlat} {\L^I}  
\nc{\git}{/\!/}   
  
\nc{\sumpl}{ { \sum_{F \in \calf_+} } }   
  
\nc{\liet}{{\bf t}}  
\nc{\nusym}{\mathcal d}  
\nc{\hfe}{h_F^{\eta} }   
\nc{\liets}{\liet^*}   
\nc{\expmf}{e^{\isq \mu_T(F)(X) } }  
\nc{\expbfj}{e^{\isq \bfj (X) } }  
\nc{\expmfs}{e^{\isq \mu_T(F)X } }  
\nc{\sumf}{\sum_{F \in \calf_+} }   
\nc{\mf}{\mu_T(F)}  
\nc{\tck}{{C^K_1}}  
\nc{\cko}{C^K_2}  
\nc{\bom}{{\bar{\omega}}}  
\nc{\fk}{F_K}  
\nc{\ft}{F_T}  
\nc{\baromega}{{\bar{\omega}}}  
  
\nc{\sss}{s}  
\nc{\wf} { { { \mathcal W}_F } }

\nc{\bfj}{ {\beta_{F,j} } }  
\nc{\nfj}{\nu_{F,j} }  
\nc{\res}{{\rm res} }

\nc{\stabo}{n_0}  
\nc{\Jac}{\Delta}  
  
\nc{\xg}{M/\!/G}  
\nc{\xtg}{\tilde{M}/\!/G}  
\nc{\txg}{\tilde{M}/\!/G}  
\nc{\mred}{\xg}  
\nc{\xred}{\xg}

\setlength{\textwidth}{6.5in}  
\setlength{\textheight}{9.1in}  
\setlength{\evensidemargin}{0in}  
\setlength{\oddsidemargin}{0in}  
\setlength{\topmargin}{-.75in}  
\setlength{\parskip}{0.3\baselineskip}

\renewcommand{\theequation}{\thesection.\arabic{equation}}  
\newcommand{\renorm}{{ \setcounter{equation}{0} }}

\title{Cohomology pairings on  
singular quotients in geometric invariant  
theory}

\author{Lisa C. Jeffrey \\  
Mathematics Department, University of Toronto\\  
Toronto, ON M5S 3G3, Canada,\thanks{The author acknowledges support  
from NSERC and the Alfred P. Sloan Foundation.}  
 \\  
Young-Hoon Kiem\\Department of Mathematics, Stanford University\\   
Stanford, CA 94305-2060, USA
\\  
Frances Kirwan \\ Balliol College, Oxford OX1 3BJ, UK
 \\and\\  
Jonathan Woolf\\Christ's College, Cambridge, CB2 3BU, UK.\thanks{Partially
supported by the Seggie-Brown   
Research Fellowship at the University of Edinburgh.}}

\date{May  2002}  
  
\maketitle

\renorm  
\section{Introduction}  
  
Let $\xg$ denote the quotient in the sense of   
Mumford's geometric invariant theory  
\cite{GIT} of a nonsingular connected complex projective  
variety $M$ by an action of a connected  
complex reductive group $G$ which is linear with respect to  
an ample line bundle $L$ on $M$.  Such quotients often  
appear as moduli spaces or as compactifications of  
moduli spaces in algebraic geometry, and their   
topology has been studied for  
many years, stimulated in particular by the  
inspiring work of Atiyah and Bott \cite{abym} in the  
early 1980s and Witten \cite{tdgr} a decade later.  
In \cite{JK1}   
(see also \cite{GK,locquant,Martin,Martin2,Martin3,tdgr})   
formulas were obtained  
for the intersection pairings  
of cohomology classes of complementary dimensions in $H^*(\xg)$  
under the assumption that   
every semistable point of $M$ is stable. In that case the quotient $\xg$  
has only orbifold singularities and its cohomology with complex coefficients  
behaves much like that of a nonsingular projective variety; in particular its  
intersection cohomology $IH^*(\xg)$ is the same as its ordinary cohomology  
$H^*(\xg)$. (Intersection cohomology is defined with   
respect to the middle perversity  
throughout this paper, and all cohomology and   
homology groups have complex coefficients.)
In this paper we shall give formulas (see Theorem 8.4) for the pairings  
of intersection cohomology classes of complementary dimensions in   
the intersection cohomology $IH^*(\xg)$ of  
geometric invariant theoretic quotients $\xg$   
for which semistability is not necessarily the same as stability  
(although we make some weaker assumptions on the action).  
We also give formulas for intersection pairings on  
resolutions of singularities  
(or more precisely partial resolutions, since orbifold  
singularities are allowed) of the quotients $\xg$ (see Proposition 8.1).  
  
Let $K$ be a maximal compact subgroup of the reductive group $G$, and let  
$\liek$ denote its Lie algebra. Then using the given linearization  
of the $G$-action on $M$ we can choose a $K$-invariant K\"{a}hler form $\omega$ on $M$,  
and the action of $K$ on $M$ is Hamiltonian with respect to the symplectic  
structure given by $\omega$; i.e. there exists a moment map $\mu:M\to  
\lieks$ for the action   
(see for example Chapter 2 of \cite{Ki1}). The inclusion of $\zloc$  
in the set $M^{ss}$ of semistable points of $M$ induces a homeomorphism  
from the Marsden-Weinstein reduction, or symplectic quotient, $\zloc /K$ to  
the geometric invariant theoretic quotient $\xg$ (\cite{Ki1} 8.14). The  
condition that semistability equals stability is equivalent to the condition  
that $0$ is  
a regular value of the moment map, and implies that the cohomology  
$H^*(\xg)$ of the quotient $\xg$ is naturally isomorphic to the equivariant  
cohomology $\hk(\zloc)$ of $\zloc$ (recall that we are working with cohomology  
with complex coefficients). The restriction map $\hk(M) \to \hk(\zloc)$ is  
surjective (\cite{Ki1} 5.4), and, as a module over  
the equivariant cohomology of a point, which we write  
as $\hk$, the equivariant cohomology $\hk(M)$ of $M$ is just the tensor  
product of its ordinary cohomology $H^*(M)$ and $\hk$ (\cite{Ki1} 5.8).   
When semistability coincides with stability, the composition of the restriction map  
$\hk(M) \to \hk(\zloc)$ and the isomorphism $\hk(\zloc) \to H^*(\xg)$ gives  
us a natural surjective ring homomorphism  
\begin{equation}   
\label{1.1} \k_M :\hk(M) \to H^*(\xg).   
\end{equation}  
  
The residue formula of \cite{JK1} is a formula for pairings of cohomology classes  
of complementary dimensions in $\xg$ in terms of equivariant cohomology  
classes in $M$ which represent them via this surjection $\k_M$,  
in the case when semistability equals stability. This formula was obtained  
from a version of Witten's nonabelian localization principle \cite{tdgr}  for  
compact Hamiltonian group actions.  
For $\epsilon >0$ and $\zeta$ a formal $K$-equivariant cohomology  
class on $M$ given by a sum $\zeta = \sum_{j \geq 0} \zeta_j$ where $\zeta_j \in  
H_K^j(M)$,  Witten defines an integral $I^{\e}(\zeta)$, which depends  
on choosing a fixed invariant inner product $<,>$ on $\liek$.  
If we represent elements of $\hk(M)$ by polynomial functions on $\liek$ with values in the   
De Rham complex $\Omega^*(M)$, this integral is  
\begin{equation} \label{1.1A} I^{\e}(\zeta) = \frac{1}{(2\pi)^s \vol(K)} \int_{X \in \liek}   
[d X] e^{-\e<X,X>/2}  
\int_M \zeta(X) \end{equation}  
where $\int_M : \hk(M) \to \hk$ is the pushforward map given by integration over  
$M$, the measure $[d X]$ on $\liek$ is induced by the fixed inner product, $\vol(K)$  
is the integral of the induced volume form on $K$ and $s$ is the dimension of $K$.  
When $\zeta$ is of the form   
$$\zeta = \eta e^{i\bom}$$  
where $\eta \in \hk(M)$ 
and $\bom$  
is the extension $\bom = \omega+ \mu$  of the symplectic form $\omega$ to an  
equivariantly closed differential form on $M$,   
Witten expresses this integral  
as a sum of local contributions. If we assume for simplicity that the stabilizer in $K$ of  
a generic point of $\mu^{-1}(0)$ is trivial, then one of these local contributions reduces to  
the evaluation  
$$\k_M(\zeta) e^{\e \Theta} [\xg]$$  
of the fundamental class $[\xg]$ against the product $\k_M(\zeta) e^{\e \Theta}$,  
where $\k_M(\zeta)$ is the cohomology class on $\xg$ induced by $\zeta$ and $\Theta$  
is the image under the natural map  
\begin{equation} \label{eq1} \hk \to \hk(\zloc) \cong H^*(\xg) \end{equation}  
of a distinguished class in $H^4_K$. In fact if we identify $\hk$ in the   
natural way with the  
space of $K$-invariant polynomials on $\liek$ then $\Theta$ is the class  
represented by $X \mapsto - \frac{1}{2}<X,X>$.   
This local contribution  
is thus a polynomial in $\e$, whereas the other terms in Witten's   
expression for $I^{\e}(\zeta)$ as  
a sum of local contributions all tend to zero exponentially fast as $\e$ tends to $0$.

Of course if the degree of $\eta \in \hk(M)$ is   
equal to the real dimension of $\xg$,  
then for $\zeta = \eta e^{i\bom}$ we have  
$$\k_M(\zeta) e^{\e \Theta}[\xg] = \k_M(\eta) e^{i\omega_0 + \e \Theta}[\xg] =   
\k_M(\eta) [\xg]$$  
where $\omega_0$ is the induced symplectic form on $\xg$;   
note that $\bom = \omega $ on $\zloc$ since $\mu$ vanishes there. Also $\k_M$ is a ring homomorphism,  
so if $\k_M(\a)$ and $\k_M(\b)$ are cohomology  
classes of complementary degrees in $\xg$ then their   
intersection pairing is given by $\k_M(\a\b)[\xg]$.  
Thus  
the behaviour as $\e \to 0$ of the integrals $I^{\e}(\a \b e^{i\bom})$   
determines the intersection pairings of cohomology classes   
$\k_M(\a)$ and $\k_M(\b)$ of  
complementary degrees in $\xg$. Moreover since $\xg$ has at worst  
orbifold singularities, its cohomology $H^*(\xg)$ satisfies Poincar\'{e}  
duality and so these pairings in principle determine the kernel of the surjection  
$\k_M: \hk(M) \to H^*(\xg),$  
and hence the ring structure of $H^*(\xg)$ given the ring structure of $\hk(M)$.  
  
In \cite{JK1} the integral $I^{\e}(\eta e^{ i\bom})$ is rewritten as an integral over the Lie  
algebra $\liet$ of a maximal torus $T$ of $K$. The localization theorem for  
compact {\em abelian} actions proved by Berline and Vergne \cite{BV1} and  
by Atiyah and Bott \cite{abmm} is used to decompose this integral as a sum of terms  
indexed by the set $\calf$ of components of the fixed point set $M^T$ of $T$ on $M$. This  
leads to a formula (the residue formula,   
Theorem 8.1 of \cite{JK1}; see Theorem 3.1 of \cite{arbrank} for a corrected  
version)   
for $\k_M(\eta) e^{i\omega_0}[\xg]$.   
In fact there is no need to include the factor of $i$ here, so we shall follow the  
conventions of \cite{arbrank} and omit it.  
If $n_0$ is the order of the  
stabilizer in $K$ of a generic point of $\zloc$ the residue formula then is  
\begin{equation} \label{res} \k_M(\eta) e^{\omega_0}[\xg] =   
\frac{n_0 (-1)^{s+n_+}}{|W| \vol(T)} \res \bigl( \cald(\xvec)^2  
\sum_{F \in \calf_+} \int_F   
\frac{i_F^*(\eta e^{\bom})(\xvec)}{e_F(\xvec)} [d\xvec] \bigr ),\end{equation}  
where $\vol(T)$ and $[d\xvec]$ are the volume of $T$ and the measure   
on its Lie algebra  
$\liet$ induced by the restriction to   
$\liet$ of the fixed inner product on $\liek$, while  
$W$ is the Weyl group of $K$, the polynomial function  
$\cald(\xvec)$ of $\xvec \in \liet$ is the product of   
the positive roots\footnote{In this paper, as in \cite{arbrank},  
we adopt the convention that weights $\beta \in \liets$ send the integer  
lattice $\L^I = \ker(\exp:\liet \to T)$ to $\ZZ$ rather than to $2\pi \ZZ$, and that  
the roots of  $K$ are the nonzero weights of its complexified adjoint  
action. This is one reason why the constant in the residue formula  
above differs from that of \cite{JK1}  
Theorem 8.1 (see the footnotes on pages 123-5 of \cite{arbrank}).} of $K$  
and $n_+ = (s-l)/2$ is the number of those positive roots; as before, $s$ is the  
dimension  of $K$ and $l$ is the dimension of $T$. Also   
$\calf_+$ is a subset of $\calf$ consisting of those components $F$ of  
the fixed point set $M^T$ on which the constant  
value taken by the $T$-moment map $\mu_T:M\to\liets$  
(which is the composition of  
$\mu:M \to \lieks$ with the natural map $\lieks \to \liets$)  lies in a certain  
cone in $\liets$ with its vertex at 0, and if  
$F \in \calf$  then $i_F:F \to M$ is the  
inclusion and $e_F$ is the equivariant Euler class of the normal bundle to $F$ in  
$M$.

In this paper we consider the more general situation where there may be  
semistable points of $M$ which are not stable (or equivalently 0 is not  
a regular value of $\mu$); we assume only that there do  
exist some stable points (or equivalently that there exist some points  
in $\zloc$ where the derivative of $\mu$ is surjective). Then there is no longer  
a natural surjection from $\hk(M)$ to $H^*(\xg)$, and since $\xg$ is in general  
singular (with singularities more serious than orbifold singularities) its  
cohomology $H^*(\xg)$ may not satisfy Poincar\'{e} duality.   
However even for singular complex projective varieties, the  
intersection cohomology groups defined by Goresky and MacPherson  
\cite{GM1,GM2} satisfy Poincar\'{e} duality, as well as the other  
properties of the cohomology groups of nonsingular complex projective  
varieties known collectively as the K\"{a}hler package. Moreover the intersection  
cohomology $IH^*(\xg)$ of the quotient $\xg$ is a direct summand of the  
ordinary cohomology of any resolution of singularities of $\xg$; this is a special  
case of the decomposition theorem of Beilinson, Bernstein, Deligne and Gabber  
\cite{BBDG}.  
  
There is a canonical  
procedure (see \cite{Ki2}) for constructing a partial resolution of singularities  
$\txg$ of the quotient $\xg$. This involves blowing $M$ up along   
a sequence of nonsingular  
$G$-invariant subvarieties, all contained in the complement $M-M^s$ of   
the set $M^s$ of stable points of $M$,  
to eventually obtain a nonsingular projective variety $\tilde{M}$ with a linear $G$-action,  
lifting the action on $M$, for which every semistable point of $\tilde{M}$ is stable.  
Then the quotient $\txg$ has only orbifold singularities, and  
the blowdown map $\pi:\tilde{M} \to M$ induces a birational morphism $\pi_G:  
\txg \to \xg$ which  
is an isomorphism over the dense open subset $M^s/G$ of $\xg$.  
  
Since we are working with complex coefficients and neglecting torsion,   
orbifold singularities cause few difficulties and in particular the  
intersection cohomology $IH^*(\xg)$ of $\xg$ is a direct summand of the  
cohomology of its partial resolution of singularities $\txg$.   
So we can consider the composition   
\begin{equation} \label{eq2} \hk(M) \to \hk(\tilde{M}) \to H^*(\txg) \to IH^*(\xg) \end{equation}  
of maps, of which the first  
is induced by the blowdown  
map $\tilde{M} \to M$, the second is $\k_{\tilde{M}}$  
(see (\ref{1.1})) and the third is the  
projection of $H^*(\txg)$ onto its direct summand $IH^*(\xg)$. This  
composition is surjective (see \cite{Ki3, Woolf2})  
and in many ways it is a natural generalization of the map  
$\k_M: \hk(M) \to H^*(\xg)$ defined when $M^{ss} = M^s$ at (\ref{1.1}), so  
we shall call it $\k_M$ too.  
Since the inclusion    
of $IH^*(\xg)$ as a direct summand of $H^*(\txg)$ respects the intersection  
pairings of classes of complementary dimensions (see \cite{Kiem}  Section  6), it is reasonable  
to hope that the residue formula (1.4)  
can be applied to the quotient $\txg$ to yield a formula for the  
intersection pairings of classes $\k_M(\a)$ and $\k_M(\b)$ of complementary  
dimensions in $IH^*(\xg)$.  
  
Unfortunately various complications arise when we try to apply the   
residue formula (\ref{res}) above to   
$\tilde{M}$ to obtain pairings on the partial desingularization $\xtg$.  
In particular, although   
the construction  
of $\txg$ from the linear $G$-action on $M$ is canonical and explicit,   
the construction  
of $\tilde{M}$ is not. In fact the procedure given in \cite{Ki2} is to blow up   
the  
set $M^{ss}$ of semistable points of $M$ along a sequence of nonsingular  
$G$-invariant closed subvarieties $V$ of $M^{ss}$, after each blow-up throwing  
out any points which are not semistable, to eventually arrive at $\tilde{M}^{ss}$  
and thus obtain $\xtg = \tilde{M}^{ss}/G$. The variety $\tilde{M}$ itself can   
be obtained by resolving the singularities of the closures $\bar{V}$ of these  
subvarieties $V$ and blowing up along their proper transforms, but in practice this  
is not usually simple. Since the residue formula (1.4) involves the set  
of components of the fixed point set of the action of the maximal torus $T$ of $K$,  
applying it directly to $\tilde{M}$ is likely to be very complicated;  
knowledge of the set of semistable points $\tilde{M}^{ss}$ alone would certainly  
not be sufficient.  Luckily it turns out that there is a simpler way to  
obtain the pairings.  
  
It is worth observing, however, that when $K$ is itself a compact torus $T$ (or equivalently when the complexification  
$G=K_c$ of $K$ is a complex torus $T_c$) then most of the difficulties described  
above disappear. In this case $\tilde{M}$ can be obtained from $M$ by  
blowing up along the components which meet $M^{ss}$  
(or equivalently which meet $\zloc$) of the fixed point  
sets of subtori $T'$ of $T$, in decreasing order of the dimension of $T'$.  
  
In the general case we can make use of the key observation due to S. Martin   
\cite{Martin,Martin2} and   
independently to Guillemin and Kalkman \cite{GK} that when  
$M^{ss} = M^s$ the evaluation $\k_M(\eta)[\xg]$ of the induced cohomology  
class $\k_M(\eta) \in H^*(\xg)$ on the fundamental class of $\xg$ is equal to  
a constant multiple, independent of $\eta\in\hk(M)$, of the evaluation on  
the fundamental class of $M/\!/T_c$ of the cohomology class  
$\k^T_M(\eta {\cald}^2 )$ on $M/\!/T_c$ induced by   
$\eta {\cald}^2 \in H^*_T(M)$.  
Here $\eta$ and $\cald$ are regarded as elements of $H^*_T(M)$ via  
the natural maps $H^*_K(M) \to H^*_T(M)$ and $H^*_T \to H^*_T(M)$. Indeed  
it follows from (\ref{res}) that if $n_0^T$ is the order of the stabilizer in $T$ of  
a generic point of $\mu_T^{-1}(0)$ then  
\begin{equation} \label{ktot} \k_M(\eta)[\xg]=\frac{n_0 (-1)^{n_+}}{n_0^T |W|}  
\k_M^T(\eta\cald^2)[M/\!/T_c] \end{equation}  
when $M^{ss} = M^s$ (although we have to be careful how we interpret  
the right hand side of this equation if semistability is not the same  
as stability for the torus  
action), and Martin \cite{Martin2} has given a direct proof of this which also shows that  
$\k_M^T(\eta\cald)[M/\!/T_c]$ is, up to a sign $(-1)^{n_+}$ which depends on a  
choice of orientation, the evaluation on the fundamental class of $\zloc/T$ of  
the cohomology class induced by $\eta \in \hT(M)$. If $M^{ss} \neq M^s$ then we can   
apply (\ref{ktot}) to the blow-up   
$\tilde{M}$ of $M$  obtained in the construction of the  
partial desingularization $\tilde{M}/\!/G$    
of $\xg$.  
  
Next we use the second stage of the approach to nonabelian localization  
due to Guillemin-Kalkman and to Martin, which involves studying the  
symplectic quotients $\mu_T^{-1}(\xi)/T$ as $\xi$ varies in $\liets$. Since  
$T$ is abelian, $\mu_t - \xi:M \to \liets$ is a moment map for the action  
of $T$ on $M$ and $\mu_T^{-1}(\xi)/T$ is a symplectic quotient which,  
when $\xi$ is rational, can be identified with the geometric invariant  
theoretic quotient $M/\!/T_c$ of $M$ by the complex torus $T_c$ with respect to a modified  
linearization.  
  
Recall that the image $\mu_T(M)$ of the moment map $\mu_T$ for the  
action of the maximal torus $T$ of $K$ is a convex  
polytope; indeed $\mu_T$ is constant on the connected components $F \in   
\calf$ of  
the fixed point set $M^T$ for the action of $T$ on $M$, and so  
$\mu_T(M^T)$ is a finite set, whose convex hull is $\mu_T(M)$   
\cite{Aconv,GSconv}.   
The convex polytope $\mu_T(M)$ is   
a union of subpolytopes, each of which is   
the convex hull of a subset of the finite set $\mu_T(M^T)$ and contains no   
points  
of $\mu_T(M^T)$ in its interior. The interior of each such subpolytope  
consists of regular values of $\mu_T$. The boundaries (or \lq walls')   
between subpolytopes  
consist of the critical values of the moment map $\mu_T$, and are the  
images under $\mu_T$ of the fixed point sets of one-parameter subgroups of  
$T$.  The evaluation on the fundamental class of  $\mu_T^{-1}(\xi)/T$ of   
the cohomology class induced by any $\eta \in \hT(M)$ is unchanged as $\xi \in \liets$   
varies within a connected component of the set of regular values of $\mu_T$,  
and in \cite{GK} formulas are obtained for the change in this evaluation  
as $\xi$ crosses a wall.   
Applying these formulas to the blow-up $\tilde{M}$ of  
$M$, we find that it is possible to choose $\xi \in \liets$ which is a regular  
value of both $\mu_T$ and $\tilde{\mu}_T$, with the following two  
properties. Firstly the difference between  
$$\k_{\tilde{M}}^T(\eta\cald^2)[\tilde{M}/\!/T_c]$$  
and   
the evaluation on the fundamental class of    
$\tilde{\mu}_T^{-1}(\xi)/T=\tilde{M}/\!/_{\xi} T_c$ of the cohomology class  induced by $\eta \cald^2 \in \hT(M)$   
can be calculated in terms of data determined purely by the construction of   
$\tilde{M}^{ss}$ from $M^{ss}$, which is canonical and explicit, rather  
than the construction of $\tilde{M}$ from $M$, which is neither canonical nor explicit.  
Secondly this evaluation on  
$[\tilde{M}/\!/_{\xi} T_c]$ is equal to the evaluation on the fundamental class of  $\mu_T^{-1}(\xi)/T$  
of the cohomology class  induced by $\eta \cald^2 $, which can be calculated by  
using the residue formula (\ref{res}) applied to the action of $T$ on $M$ with  
the moment map $\mu_T - \xi$. This combined with (\ref{ktot}) enables us to   
calculate pairings in the cohomology of the partial desingularization  
$\xtg$ of $M/\!/G$.  
  
In order to understand pairings in $IH^*(\xg)$ of intersection cohomology  
classes on the singular quotient $\xg$ we make use of the work of  
the second author in \cite{Kiem}.   
First note that the composition $\k_M:\hk(M) \to  
IH^*(\xg)$ at (\ref{eq2}) factors as the composition  
of the restriction map $\hk(M) \to \hk(M^{ss})$ and  
a surjection   
\begin{equation} \label{kss}  
\k^{ss}_M : \hk(M^{ss}) \to IH^*(\xg).\end{equation}  
In \cite{Kiem} it is  
shown that if the action of $G$ on $M$  
is weakly balanced (see Definition 5.3 below), then there is a naturally  
defined subset $V_M$ of $\hk(M^{ss})$ such that  
$\k_M^{ss} : \hk(M^{ss}) \to IH^*(\xg)$ restricts to an isomorphism  
$$\k_M^{ss}: V_M \to IH^*(\xg).$$  
In \cite{Kiem} it is also shown that the intersection pairing of two elements  
$\k_M(\alpha)$ and $\k_M(\beta)$ of complementary degrees in $IH^*(\xg)$  
is equal to the evaluation of the image in $H^*(\tilde{M}/\!/G)$ of the  
product $\alpha\beta \in \hk(M)$ on the fundamental class  
$[\tilde{M}/\!/G]$, provided that $\alpha|_{M^{ss}}$ and $\beta|_{M^{ss}}$ lie in $V_M$.  
In Section 8 below we shall show that if $\alpha|_{M^{ss}}$ and $\beta|_{M^{ss}}$ lie in $V_M$  
and are of complementary degrees with respect to $\xg$, then the intersection  
pairing of $\k_M(\alpha)$ and $\k_M(\beta)$ in $IH^*(\xg)$ is given, just as at (\ref{res}),  
by  
\begin{equation}  
\langle \k_M(\a),\k_M(\b)\rangle_{IH^*(\xg)} =   
\frac{n_0 (-1)^{s+n_+}}{|W| \vol(T)} \res \bigl( \cald(\xvec)^2  
\sum_{F \in \calf_+} \int_F \frac{i_F^*(\a\b e^{\bom})  
(\xvec)}{e_F(\xvec)} [d\xvec] \bigr),  
\end{equation}  
provided that the multivariable residue $\res$ and  
the subset $\calf_+$ of $\calf$ are   
interpreted correctly. In the case when $T$ is  
one-dimensional, as before we can take $\calf_+$ to be the  
set of those $F \in \calf$ such that $\mu_T(F)$ is positive. The difference   
now is that there may be some $F \in \calf$ with   
$\mu_T(F)=0$, which cannot happen when semistability  
coincides with stability; this suggests that we need to  
be careful to decide whether $\calf_+$ consists  
of those $F\in \calf$ for which $\mu_T(F)$ is   
non-negative, or just those for which $\mu_T(F)$   
is strictly positive. However it turns out that  
when $\a$ and $\b$ lie in $V_M$ then  
$$\res \bigl( \cald(X)^2 \sum_{F \in \calf, \,\,  
\mu_T(F)=0} \int_F \frac{i_F^*(\a\b e^{\bom})(\xvec)}  
{e_F(\xvec)} [d\xvec] \bigr) = 0,$$  
so in fact it does not matter which definition of  
$\calf_+$ we choose here, and the situation is similar   
when $\dim T > 1$.  
  
We can also consider Witten's integral $I^{\e}(\eta e^{i\bom})$. When 0  
is a regular value of the moment map $\mu$ (or equivalently when  
semistability is the same as stability) then, as we have seen,  
$I^{\e}(\eta e^{i\bom})$ can be expressed as a sum of terms which tend to 0   
exponentially fast with $\e$, together with   
$$ \k_M(\eta) e^{i\omega_0+\e\T}[\xg] $$  
which is a polynomial in $\e$ and can be expressed  
using the residue formula (1.4) as a sum over the components  
$F \in \calf$ of $M^T$. When 0 is not a regular value of $\mu$ we can  
still write Witten's integral $I^{\e}(\eta e^{i\bom})$  as a sum  
of exponentially small terms together with a sum over the  
components $F \in \calf$ of $M^T$ (see Section 9 below). The terms in this sum indexed  
by $F\in \calf$ such that $\mu_T(F)$ does not lie on a wall through 0  
are exactly as they would be in the case when 0 is a regular value of $\mu$, i.e. the  
residue of  
$$\frac{n_0 (-1)^{s+n_+}}{|W| \vol(T)}  \bigl( \cald(\xvec)^2  
e^{-\e <\xvec,\xvec>/2} \int_F   
\frac{i_F^*(\eta e^{i\bom})(\xvec)}{e_F(\xvec)} [d\xvec] \bigr ). $$  
In particular these terms are polynomials in $\e$. However the terms  
indexed by $F\in \calf$ such that $\mu_T(F)$ does lie on a wall through 0  
are in general  
only polynomials in $\e^{1/2}$, and it is unclear whether the  
sum can be interpreted in terms of intersection pairings when  
0 is not a regular value of $\mu$.

The construction of the partial desingularization $\txg$ can  
also be carried out in the symplectic category, using symplectic  
blow-ups, to give a partial desingularization of the symplectic  
reduction of any Hamiltonian $K$-action on a compact symplectic  
manifold \cite{MSj,Woolf}. Because symplectic blow-ups depend  
on a number of choices the partial desingularizations obtained  
will not be unique up to symplectomorphism, but they will be  
determined up to symplectic homotopy, and in particular up to  
diffeomorphism. The analysis of Witten's  
integral $I^{\e}(\eta\e^{i\bom})$ and the formulas for pairings in $H^*(\xtg)$ and $IH^*(\xg)$ are  
also valid for singular  
symplectic reductions. 
  
In \cite{tdgr} Witten studied the moduli spaces $\calm (n,d)$ of holomorphic bundles  
of rank $n$ and degree $d$  
over a fixed compact Riemann surface $\Sigma$ as symplectic reductions  
of infinite dimensional affine spaces by infinite dimensional Lie  
groups. When the rank $n$ and degree $d$ of the bundles are coprime  
(i.e. when semistability is the same as stability  
and the moduli space $\calm(n,d)$ is nonsingular) then   
using physical methods Witten  
obtained formulas (later proved using different methods in  
\cite{arbrank}) for intersection pairings on these moduli spaces  
from asymptotic expansions of the integrals  
$I^{\e}(\eta e^{i\bom})$ as $\e$ tends to 0. He also gave formulas  
for the asymptotic expansions of the integrals in the simplest case when semistability differs from stability, namely the case of bundles  
of rank two and even degree, and he noted that powers of $\e^{1/2}$   
appeared. In a forthcoming article \cite{JKKW} we will use the finite dimensional  
methods of \cite{arbrank} together with the results of this paper  
to rederive Witten's calculations   
for $\calm(2,0)$ and give formulas for intersection  
pairings in $IH^*(\calm (n,d))$ and on the partial resolution of singularities  
of the moduli space $\calm (n,d)$, in the general case for $n\geq 2$   
when $n$ and $d$ may have common factors so that  
$\calm(n,d)$ may have singularities.  
  
The layout of this paper is as follows. In Section 2 we recall briefly the   
relationship between geometric invariant theory (GIT) and the moment map  
in symplectic geometry, and the use of equivariant cohomology  
to study the cohomology of GIT quotients. In Section 3 we review   
Witten's principle of nonabelian localization and the  
residue formula of \cite{JK1} in the case when 0 is a regular  
value of the moment map. In Section 4 we recall the construction of the  
partial desingularization $\txg$. In Section 5 we review intersection   
cohomology and the work of the second author in \cite{Kiem}, and in Section 6 we  
study intersection pairings in $IH^*(\xg)$ via the  
isomorphism $\k_M^{ss}:V_M \to IH^*(\xg)$ from   
\cite{Kiem}. In Section 7  we give formulas for   
pairings in the intersection  
cohomology $IH^*(\xg)$ of the singular quotient, and in Section 8 we   
calculate pairings in the cohomology $H^*(\txg)$ of the partial desingularization.  
Finally in Section 9 we study Witten's integral.

\renorm  
\section{The moment map and cohomology of quotients}  
  
In this section we shall recall briefly the relationship between   
geometric invariant theory  
and the moment map in symplectic geometry (see e.g. \cite{GIT} or \cite{Ki1} for  
more details), and the use of equivariant cohomology to study the  
cohomology of geometric invariant theoretic quotients.  
  
Let $M$ be a nonsingular connected complex projective variety, and let $G$  
be a connected complex reductive group acting on $M$. In order  
to define the geometric invariant theoretic quotient $\xg$ we need a linearization of the action  
of $G$ on $M$; i.e. we need a lift of the action to a linear action on a line  
bundle $L$ over $M$, which is usually assumed to be ample. We shall suppose for simplicity that $M$ is embedded in  
a complex projective space $\PP_n$ and that $L$ is the hyperplane line  
bundle on $M$; then we need the action of $G$ to be given by a  
representation $\rho:G \to GL(n+1)$. The quotient $\xg$ is the projective variety  
defined by the finitely generated graded subalgebra of  
$\bigoplus_{k\geq 0} H^0(M,L^{\otimes k})$ consisting of all  
elements invariant under the action of $G$.  
  
There is a surjective $G$-invariant morphism $\tau: M^{ss} \to \xg$ from an open $G$-invariant  
subset $M^{ss}$ of $M$ (whose elements are called semistable points of $M$) to  
$\xg$; in fact $x \in M$ is semistable if and only if there is a $G$-invariant  
section of $L^{\otimes k}$ for some $k$ which does not vanish at $x$. If $x$ and $y$ are semistable points of $M$ then $\tau(x) = \tau(y)$ if and only  
if the closures of the orbits $Gx$ and $Gy$ meet in $M^{ss}$. There is an open   
$G$-invariant subset $M^s$ of $M^{ss}$ (whose elements are called stable\footnote{This  
is now the usual terminology and notation. However in \cite{GIT} the terminology  
\lq\lq properly stable'' and notation $M^s_{(0)}$ are used instead.} points of $M$) such  
that every fibre of $\tau$ which meets $M^s$ is a single $G$-orbit of dimension equal  
to the dimension of $G$. We shall assume that $M^s$ is nonempty; then the image  
of $M^s$ in $\xg$ is open and dense and can be identified naturally with the  
ordinary topological quotient $M^s/G$.  
  
We shall call elements of $M^{ss} - M^s$ {\it strictly semistable}, and  
write $M^{sss} = M^{ss} - M^s$.  
  
The subsets $M^{ss}$ and $M^s$ of $M$ are characterized by the following properties (see  
Chapter 2 of \cite{GIT} or \cite{New}).   
  
\begin{prop} \label{sss} (i) A point $x \in M$ is semistable (respectively  
stable) for the action of $G$ on $M$ if and only if for every  
$g\in G$ the point $gx$ is semistable (respectively  
stable) for the action of a fixed maximal (complex) torus of $G$.  
  
\noindent (ii) A point $x \in M$ with homogeneous coordinates $(x_0,\ldots,x_n)$  
in some coordinate system on $\PP_n$  
is semistable (respectively stable) for the action of a maximal (complex)  
torus of $G$ acting diagonally on $\PP_n$ with  
weights $\a_0, \ldots, \a_n$ if and only if the convex hull  
$$\conv \{\a_i :x_i \neq 0\}$$  
contains $0$ (respectively contains $0$ in its interior).  
\end{prop}  
  
Now let $K$ be a maximal compact subgroup of $G$; then $G$ is the  
complexification of $K$. By choosing coordinates on $\PP_n$ appropriately  
we may assume that $K$ acts unitarily. Then $K$ preserves the K\"{a}hler  
structure on $M$ given by the restriction of the Fubini-Study metric  
on $\PP_n$. The K\"{a}hler form $\omega$ makes $M$ into a symplectic  
manifold on which $K$ acts. Associated to this action there is a moment  
map  
$\mu:M\to \lieks,$  
where $\liek$ is the Lie algebra of $K$, given in homogeneous  
coordinates $x=(x_0,\ldots,x_n)$ by the formula  
$$\mu(x).a = (2\pi i |\!| x |\!|^2 )^{-1} x \rho_*(a) \bar{x}^t $$  
for $a \in \liek$. Here the element $\rho_*(a)$ of the Lie algebra  
of the unitary group $U(n+1)$ is thought of as an $n+1$ by $n+1$  
skew-hermitian matrix. Recall that the defining property of a moment map $\mu:M\to  
\lieks$ is that  
\begin{equation} \label{2.1} d\mu(x)(\xi).a = \omega_x(\xi, \tilde{a}(x)) \end{equation}  
for all $x \in M$, $\xi \in T_xM$ and $a\in\liek$, where $\tilde{a}$ is the  
vector field on $M$ induced by $a$. We also require that $\mu$ carries the  
given $K$-action on $M$ to the coadjoint action on the dual of its Lie algebra.

When a compact group $K$ acts on a symplectic manifold $M$ and  
$\mu:M \to \lieks$ is a moment map for the action, the  
symplectic form on $M$ induces a symplectic form  
on the quotient $\zloc/K$   
(away from its singularities, at least),  
which  
is the Marsden-Weinstein reduction or  
symplectic quotient of $M$ by the action of $K$. In our situation $\zloc/K$ can be  
identified with the geometric invariant theoretic quotient $\xg$.  
A more precise statement is the following (see \cite{GIT} Theorem 8.2 and Remark 8.3   
or \cite{Ki1} 6.2, 8.10, 7.2 and 7.5).  
  
\begin{proposition} \label{2.2} (i) $x \in M^{ss}$ if and only if $\zloc$ meets the   
closure of the orbit $Gx$ in $M$.  
  
\noindent (ii) $x \in M^s$ if and only if $\zloc$ meets  
the orbit $Gx$ at a point whose stabilizer in $K$ is finite.  
  
\noindent (iii) The inclusion of $\zloc$ in $M^{ss}$ induces a homeomorphism  
$\zloc/K \to \xg$.  
\end{proposition}  
  
If the Lie algebra $\liek$ is given a fixed $K$-invariant inner product then we can consider  
the function $|\!| \mu|\!|^2$ as a Morse function on $M$ (although it is not  
a Morse function in the classical sense; see \cite{Ki1}). It induces a Morse stratification  
of $M$, in which the stratum containing any $x \in M$ is determined by the limit  
set of its path of steepest descent for $|\!| \mu|\!|^2$ (with respect to the K\"{a}hler  
metric). This stratification can also be defined purely algebraically, and has the following  
properties (see \cite{Ki1} 5.4 and Chapters 12 and 13).  
  
\begin{proposition} \label{2.3} (i) Each stratum is a $G$-invariant locally  
closed nonsingular subvariety of $M$.  
  
\noindent (ii) $M^{ss}$ coincides with the unique open stratum.  
  
\noindent (iii) The stratification is equivariantly perfect   
(that is, its equivariant Morse inequalities are all equalities)  
over the complex numbers,  
and in particular the restriction map  
$$\hk(M) \to \hk(M^{ss}) \cong \hk(\zloc)$$  
is surjective.  
\end{proposition}  
  
Here the $K$-equivariant cohomology of  
any topological space $Y$ on which $K$ acts is  
$$\hk(Y) = H^*(Y \times_K EK)$$  
where $EK \to BK$ is the universal $K$-bundle. (Recall that all cohomology  
groups have complex coefficients throughout this paper). Note that $K$ is  
homotopically equivalent to its complexification $G$, so $G$-equivariant cohomology  
is naturally isomorphic to $K$-equivariant cohomology; we shall work with the  
latter in this paper.  
  
If $M$ is a manifold, the $K$-equivariant cohomology of $M$   
can be identified with the cohomology  
of a chain complex $\Omega^*_K(M)$ whose elements  are $K$-equivariant  
polynomial functions on the Lie algebra $\liek$ of $K$ with values in the  
de Rham complex $\Omega^*(M)$ of differential forms on $M$ (see for example  
Chapter 7 of  \cite{BGV}). We shall call elements of $\Omega^*_K(M)$  
equivariant differential forms on $M$. The differential $D$ on this complex  
is defined by\footnote{This  
definition of the equivariant cohomology differential differs by a factor of $i$ from  
that used in \cite{tdgr} but is consistent with that used in \cite{JK1,arbrank}.}  
\begin{equation} \label{1.0002}  
(D\eta)(\xvec) = d(\eta (\xvec) )  - \iota_{\xvec^{\#}} (\eta(\xvec) ) \end{equation}  
where $\xvec^{\#}$ is the vector field on $M$ generated by the action of   
$\xvec$ (see Chapter 7 of \cite{BGV}). We can write  
$\Omega^*_K(M)=(S(\lieks)\otimes \Omega^*(M))^K$ where $S(\lieks)$ denotes the algebra  
of  polynomial functions on the Lie algebra $\liek$ of $K$. An element   
$\eta \in \Omega^*_K(M)$ may be thought of   
as a $K$-equivariant polynomial function from $\liek$   
to $\Omega^*(M)$, or alternatively as a family of differential forms  
on $M$ parametrized by $\xvec \in \liek$.   
The standard definition of degree is used on $\Omega^*(M)$ and  
degree two is assigned to elements of $\lieks$.

In fact, as a vector space though not in general as a ring,  
when $M$ is a compact symplectic manifold with a Hamiltonian action   
of $K$ then $H^*_K(M)$ is   
isomorphic to $H^*(M) \otimes H^*_K$ where $H^*_K=\Omega^*_K({\rm pt})  
=S(\lieks)^K$ is the  
equivariant cohomology of a point (see \cite{Ki1} Proposition 5.8).  
  
It follows directly from the defining property of a moment   
map that if $\mu$  
is regarded in the obvious way as a linear map from $\liek$ to the space  
$\Omega^0(Y)$ of smooth complex-valued functions on $Y$, then  
$\bom \in \Omega^2_K(M)$ defined by  
$$\bom(X) = \omega + \mu(X)$$  
satisfies $D\bom =0$ and therefore defines an extension of the  
cohomology class of $\omega$ in $H^2(M)$ to an equivariant cohomology  
class in $H^2_K(M)$.   
  
If every semistable point of  $M$ is stable then by Proposition \ref{2.2}  
$K$ acts on $\zloc$ with only finite stabilizers. Because of the defining  
property (\ref{2.1}) of a moment map, this implies that $0$ is a regular  
value of $\mu$ and hence that $\zloc$ is a submanifold of $M$. Since  
the cohomology with complex coefficients of a classifying space of  
a finite group is always trivial, it also implies that the obvious map  
$$\zloc \times_K EK \to \zloc/K$$  
induces an isomorphism  
\begin{equation} \label{2.4} H^*(\zloc/K) \cong \hk(\zloc) \end{equation}  
and hence  
$H^*(\xg) \cong \hk(M^{ss}).$  
Composing this with the surjection of Proposition \ref{2.3}(iii), we find  
that if $M^{ss}=M^s$ then there is a natural surjective ring homomorphism  
from $\hk(M)$ to $H^*(\xg)$.

\begin{example} Consider the action of $G=SL(2)$   
and its maximal compact subgroup $K=SU(2)$ on   
$\PP_n$ identified with the space of unordered   
sequences of $n$ points in $\PP_1$ (that is,  
with the projectivized symmetric product $\PP(S^n(\CC^2))$).  
The diagonal subgroup $\CC^*$ is a maximal torus of $G$ and   
acts with weights $n,n-2,n-4,...,2-n,-n$ on $S^n(\CC^2)=\CC^{n+1}$.  
An element $[a_0,...,a_n]$ of $\PP_n$ corresponds to  
the $n$ roots in $\PP_1$ of the polynomial with  
coefficients $a_0,...,a_n$; it is semistable (respectively   
stable) for the action of $G$ if and  
only if at most $n/2$ (respectively strictly fewer  
than $n/2$) of these roots coincide  
anywhere on $\PP_1$.   
The induced stratification of $M$ has strata  
$S_0=M^{ss}$ and $S_j$ for $n/2 < j \leq n$. If $n/2 < j \leq n$  
then the elements of $S_j$ correspond to sequences of $n$ points on  
$\PP_1$ such that exactly $j$ of these points coincide somewhere on  
$\PP_1$, and $S_j$ retracts equivariantly onto  
the subset of $M$ where $j$ points coincide somewhere on $\PP_1$ and  
the remaining $n-j$ points coincide somewhere else on $\PP_1$.   
This subset is a single $G$-orbit with stabilizer $\CC^*$, so  
that   
$$H_K^*(S_j) \cong H^*(B\CC^*) \cong H^*(BS^1),$$  
and the fact that the stratification is  
equivariantly perfect tells us that  
$$\dim H^q_K(M^{ss}) = \dim H_K^q(M) - \sum_{n/2 < j \leq n} \dim H_K^{q  
- \dim_{\RR}(S_j)}(S_j).$$  
The same is true when $M=(\PP_1)^n$, except that then $S_j$ has $(^n_j)$   
components, each of  
which retracts onto a single $G$-orbit and has equivariant cohomology  
isomorphic to $H^*(BS^1)$ (see \cite{Ki1} Section 9 for more details).  
\end{example}  
  
\begin{example}\label{5.2} (Example 6.3 in \cite{Kiem}.)   
Consider the $\CC^*$-action on $\PP^7$ by a representation with  
weights $+1, 0, -1$ with multiplicity $3,2,3$ respectively. Then  
$$H^*_{S^1}(\PP^7)=\CC[\xi,\rho]/<\xi^2(\xi-\rho)^3(\xi+\rho)^3>$$  
is the quotient of the polynomial ring $\CC[\xi,\rho]$  
where $\xi$ is a generator in $H^2(\PP^7)$ and $\rho$ is a generator  
in $H^2(BS^1)$, by the ideal generated by  
$\xi^2 (\xi - \rho)^3(\xi + \rho)^3$.  
There are two unstable strata whose equivariant cohomology  
classes are $\xi^2(\xi-\rho)^3$ and $\xi^2(\xi+\rho)^3$. Since  
the Morse stratification with respect to the norm  
square of the moment map is equivariantly perfect,  
$$H^*_{S^1}((\PP^7)^{ss})  
=\CC[\xi,\rho]/<\xi^2(\xi-\rho)^3,\, \xi^2(\xi+\rho)^3>$$  
A Gr\"obner basis for the relation ideal is   
$$\{\xi^5+3\xi^3\rho^2, \xi^4\rho+\frac13\xi^2\rho^3, \xi^3\rho^3,   
\xi^2\rho^5\}$$  
where $\xi>\rho$. Hence as a vector space,  
$$H^*_{S^1}((\PP^7)^{ss})  
\cong \CC\{\xi^i\rho^j : i=0,1,\, j\ge 0\}\oplus  
\CC\{\xi^i\rho^j : 2i+j<9, i\ge 2, j\ge 0\}.$$  
\end{example}

\renorm 
\section{Residue formulas and nonabelian localization}

The map $\Omega^*_K(M) \to 
\Omega^*_K({\rm pt}) = S(\lieks)^K$ 
given by integration over $M$ passes to $\hk(M)$. Thus  for any 
$D$-closed element
$\eta \in \Omega^*_K(M)$ representing a cohomology class $[\eta]$, 
there is a corresponding element 
$\int_M \eta \in \Omega^*_K({\rm pt})$ which depends only on 
$[\eta]$. The same is true for any 
$D$-closed formal series  $\eta = \sum_j \eta_j$ 
of elements $\eta_j$ in $\Omega^j_K(M)$: we shall in particular consider terms
of the form 
$$\eta (\xvec) e^{ \bom(\xvec)} $$ 
where 
$\eta \in \Omega^*_K(M)$ and 
$$\bom(\xvec) = \omega + \mu(\xvec) \in \Omega^2_K(M). $$

 If $\xvec$ lies in $\liet$, the Lie algebra of 
a maximal torus $T$ of $K$, then there is a formula for 
$\int_M \eta(\xvec)$ (the {\em   abelian localization 
formula} \cite{abmm,BGV,BV1,BV2})
which depends only on the fixed point set of 
$T$ in $M$. It tells us that 
\begin{equation} \label{1.002}
 \int_M \eta(\xvec) 
= \sum_{F \in \calf} \int \frac{i_F^* \eta(\xvec)}
{e_F(\xvec)} \end{equation}
where $\calf$ indexes the components $F$ of the fixed point set of 
$T$ in $M$, the inclusion of $F$ in $M$ is denoted
by $i_F$ and $e_F $ 
$\in H^*_T(M)$ is the equivariant Euler class of the normal 
bundle to $F$ in $M$. In particular, applying (\ref{1.002})
with  $\eta $ replaced by the formal equivariant cohomology class
 $\eta e^{ \bom}$ 
we have 
\begin{equation} \label{1.003} 
 \int_M \eta(\xvec)e^{ \bom(\xvec)}  
= \sum_{F \in \calf} h^\eta_F(\xvec), \end{equation}
where 
\begin{equation} \label{1.004}   \hfeta(\xvec) = 
e^{ \mu(F)(\xvec)}\int_F  \frac{i_F^* \eta(\xvec)  e^{ \omega} 
 }{e_F (\xvec) }. \end{equation}
Note that the moment map $\mu$ takes a 
constant value $\mu(F) = \mu_T(F) \in\liets$ 
 on each 
$F \in \calf$, and that the integral in (\ref{1.004}) 
is a rational function of $\xvec$. 

The main result (the residue formula, Theorem 
8.1) of \cite{JK1} gives a formula for the evaluation on the fundamental
class $[\mred]\in H_*(\mred)$, or equivalently
(if we represent cohomology classes by differential forms) the integral
over $\mred$, of the image
$\k_M(\eta) e^{\omega_0}$ in $H^*(\mred)$ of any
formal equivariant cohomology class on $M$ of the type $\eta e^{\bom}$ where
$\eta \in \hk(M)$. 

\begin{theorem} \label{t4.1}{\bf  (i) (Residue 
formula, \cite{JK1} Theorem 8.1)}  Let $\eta \in 
\hk(M) $ induce $\k_M(\eta) \in H^*(\xred)$. Then 
\begin{equation} \label{jk81}   \k_M(\eta) e^{{}\omega_0} [\xred] 
 = n_0 {C_K}  \res \Biggl ( 
\cald^2 (X)
 \sum_{F \in \calf} \hfeta(X) [d X] \Biggr ), \end{equation}
where the constant\footnote{In this paper we are adopting the conventions of
\cite{arbrank} on the equivariant differential and on the
normalization of weights 
(see Footnote 9 on page 125 of \cite{arbrank} for the
effect of different conventions on the constant $C_K$).}
$C_K$ is defined by 
\begin{equation} \label{4.001} C_K = \frac{(-1)^{s+n_+}}{ |W| \vol(T)}, \end{equation}
and $n_0$ is the order of the stabilizer in $K$ of a generic 
point of $\zloc$.

{\bf (ii) (Reduction to an integral over $\mu^{-1}(0)/T$, 
\cite{Martin})}
$$\k_M(\eta) e^{{}\omega_0} [\xred] = 
\frac{n_0 }{n_0^T |W| } \int_{\mu^{-1} (0)/T} 
\kappa_M^T (  
\cald (X) \eta e^{{} \bar{\omega} } ),$$
where $n_0^T$ is the order of the stabilizer in 
$T$ of a generic point of $\mu^{-1}(0)$. 
\end{theorem}

In these formulas $|W|$ is the order of the Weyl group $W$ of $K$, 
while $s = \dim K$ and 
$l  = \dim T$, and
$n_+ = (s-l)/2$ is the 
number of positive roots. The measure $[dX]$ on $\liet$ and volume $\vol(T)$ of
$T$ are obtained from the restriction of a fixed invariant inner product on $\liek$,
which is used to identify $\lieks$ with $\liek$ throughout.
Also, $\calf$ denotes the set of components of the fixed point
set of $T$, and if $F$ is one 
of these components then the meromorphic  function 
$\hfeta$ on $\liet \otimes \CC$ is defined by 
(\ref{1.004}). The polynomial
$\cald: \liet \to \RR$ is defined by 
$$\cald(X) = \prod_{\g > 0 } \g(X),$$
where $\g$ runs 
over the positive roots of $K$. Note that it would perhaps be more
natural to combine $(-1)^{n_+}$ from the constant $C_K$ with
$\cald^2(X)$ and replace them by the product
$$\prod_{\g} \g(X)$$
of all the positive and negative roots of $K$.

Let $\mu_T:M\to\liets$ be the composition of the moment map
$\mu:M \to \lieks$ with the restriction map from $\lieks$
to $\liets$; then $\mu_T$ is a moment map for the action of
$T$ on $M$. In particular $\mu_T$ is constant on any connected component $F$ of the fixed point set $M^T$
for the action of $T$ on $M$.

The multivariable residue $\res$ which appears in the formula (3.4) above can
be thought of as a linear map defined on a certain class of  meromorphic differential forms on $\liet\otimes\CC$,
but in order to apply it to the individual terms in the residue
formula it is necessary to make some choices which do not affect the residue
of the whole sum. Once the choices have been made, many
of the terms in the sum contribute zero and the formula
 can be rewritten as
a sum over a subset $\calf_+$ of the set $\calf $ of components of the
fixed point set $M^T$, consisting of those $F \in \calf$ on which the constant
value taken by $\mu_T$ lies in a certain cone with its vertex at 0. 
When the rank of $K$ is one and $\liet$ is identified
with $\RR$, we can take
$$\calf_+ = \{ F \in \calf: \mu_T(F) > 0 \}. $$ 
When $K = U(1)$, then the residue formula becomes
\begin{equation} \label{u1}
\k_M(\eta) e^{{}\omega_0} [\xred] = - n_0
 \res_{X=0} \Bigl (  \sum_{F \in \calf_+} \hfeta(X) dX
\Bigr ) \end{equation} 
where $\res_{X=0} $ denotes the coefficient of $1/X$ when
$X \in \RR$ has been identified with $2 \pi i X \in \liek$.  
When $K = SU(2)$ we have 
\begin{equation} \label{su2} \k_M(\eta) e^{{}\omega_0} [\xred] = 
\frac{n_0}{2} \res_{X=0} \Bigl ( (2 X)^2  \sum_{F \in \calf_+} \hfeta(X)
 \Bigr )  \end{equation} 
when $X \in \RR $ has been identified with
${\rm diag} (2 \pi i, - 2 \pi i) X \in \liet$. 

\begin{example} When $K=SU(2)$ with maximal
torus $T=S^1$ acts on $M=(\PP_1)^n$, the
equivariant cohomology
$H^*_T(M)$ of $M$ with respect to $T$ is generated by $n$ elements 
$\xi_1,...,\xi_n$ of degree two which are lifts of the
standard generators of $H^*(M)$,
together with
another generator $\zeta$ of degree two coming
from $H^*_T$, subject to the relations
$$(\xi_1)^2 = ... = (\xi_n)^2 = \zeta^2;$$
$\hk$ is generated by $\xi_1$,...,$\xi_n$
and $\zeta^2$ subject to the same relations.
We assume that $n$ is odd, so that 
0 is a regular value of $\mu$, or equivalently 
semistability coincides with stability for the
action of the complexification $G=SL(2)$ of $K$
(cf. Example 2.4). When $\PP_1$ is identified with the unit sphere
$S^2$ in $\RR^3$ and the dual of the Lie algebra
of $SU(2)$ is identified suitably with $\RR^3$ the moment
map is given by
$$\mu(x_1,...,x_n) = x_1 + ...+x_n.$$
The fixed point sets of the action of the standard maximal torus $T$ of $K$ on $M$ are the
$n$-tuples $(x_1,...,x_n) \in (\PP_1)^n$ such
that each $x_j$ is either 0 or $\infty$.
If we index these by sequences $(\delta_1,...,\delta_n)$
where $\delta_j=1$ if $x_j=0$ and $\delta_j=-1$ if $x_j=\infty$,
then (\ref{su2}) gives us the formula
$$\k_M(q(\xi_1,...,\xi_n,\zeta^2))e^{\omega_0}[\xg] = $$ 
$$
\res_{X=0} \Bigl( 4X^2 \sum_{(\delta_1,...,\delta_n)\in \{1,-1\}^n,
\delta_1 + ... + \delta_n >0} \frac{q(\delta_1 X,...,\delta_n X,
X^2)e^{(\delta_1  + ... + \delta_n )X}}{(\prod_j \delta_j)X^n}
\Bigr)$$
for any polynomial $q(\xi_1,...,\xi_n,\zeta^2)$ in
the generators $\xi_1,...,\xi_n$ and $\zeta^2$ for
$\hk(M)$ (see Section 9 of \cite{JK1}).
\end{example}

\begin{rem} If we suppose that the degree of
$\eta$ is equal to the real dimension of $\xg$ then
of course we have
$\k_M(\eta)e^{\omega_0}[\xg] = \k_M(\eta)[\xg].$
In order that the multivariable residue $\res$ which
appears in the general version (3.4) of the residue formula
should be defined, we still need to include the terms 
$e^{\mu(F)(X)}$ and $e^{\omega}$ coming from
$e^{\bom}$ in the right hand side of (3.4). However
we can omit them from the right hand side of (\ref{u1})
and (\ref{su2}) if we wish; they have done their job
in reducing the sum over $F \in \calf$ to a sum over
$F \in \calf_+$.
\end{rem}

Note that (see \cite{JK1} 2.7) the
reciprocal of the $T$-equivariant Euler class $e_F(X)$
can be expressed in the form
$$\frac{1}{e_F(X)} = \prod_{j=1}^{N_F} \frac{1}{c_1(\nu_{F,j})+ \beta_{F,j}(X)}
= \prod_{j=1}^{N_F} \sum_{r_j\geq 0} 
\frac{(-c_1(\nu_{F,j}))^{r_j}}{\beta_{F,j}(X)^{r_j+1}}$$
where $\beta_{F,1}, \ldots,\beta_{F,N_F}$ are the weights
of the action of $T$ on the normal bundle to $F$ in $M$,
and $c_1(\nu_{F,1}), \ldots,c_1(\nu_{F,N_F}) $ $\in H^2(F)$
are nilpotent. Thus the terms $\cald^2(X)\hfe(X)$ appearing
in the residue formula can all be expressed as finite
sums of functions of the form
\begin{equation} \label{f1} h(X)=\frac{q(X)e^{\l(X)}}{\prod_{j=1}^N
\beta_j(X)}, \end{equation}
where $q(X)$ is a polynomial in $X\in \liet$ and $\l(X)$ and 
$\b_1(X),\ldots,\b_N(X)$ are linear functions of $X$.
It is shown
in Proposition 3.2 of \cite{locquant} and Proposition 8.11 of \cite{JK1} that the multivariable
residue of $h(X)[dX]$ when $h(X)$
has this form is determined completely
by a few elementary properties. Alternatively $\res$ can be expressed
in terms of iterated
one-dimensional residues using Proposition 3.4 of \cite{locquant}.

\begin{rem}
Note that the multivariable residue defined and used in \cite{JK1,locquant} is very
slightly different from the one used here and in \cite{arbrank}, because in \cite{JK1, locquant}
the residue formula is applied to formal equivariant cohomology classes of the form
$\eta e^{i \bom}$ instead of $\eta e^{ \bom}$. The factors of $i$ were omitted in \cite{arbrank} 
because they are essentially irrelevant to the residue formula, although they appear naturally
in Witten's integral $I^{\epsilon}(\eta e^{ i \bom})$. In \cite{JK1, locquant}, and also Section 9
of this paper, functions of the form $q(X)e^{i \l(X)}/ \prod_{j=1}^N
\beta_j(X)$ replace the functions of the form $q(X)e^{\l(X)}/ \prod_{j=1}^N
\beta_j(X)$ studied here and in \cite{arbrank}. To obtain the elementary properties which
uniquely determine the multivariable residue used here and in \cite{arbrank}, one simply
omits all the occurrences of $i$ in \cite{locquant} Proposition 3.2 (see also Section 9 below).

\end{rem}

\renorm
\section{Partial resolution of singularities}
\label{partial resolution}

In this section we shall describe the construction of the partial resolution
of singularities $\xtg \to \xg$ (for more details see \cite{Ki2}).

As before let $M$ be a nonsingular complex projective variety embedded
in a projective space $\PP_n$ and let $G$ be a connected
complex reductive group acting on
$M$ via a representation $\rho:G \to GL(n+1)$.
Let $V$ be any nonsingular $G$-invariant closed subvariety of $M$ and let
$\pi:\hat{M} \to M$ be the blowup of $M$ along $V$. The linear action of $G$ on the
hyperplane line bundle $L$ over $M$ lifts to a linear action on the line bundle
over $\hat{M}$ which is the pullback of $L^{\otimes k}$ 
tensored with ${\mathcal O}(-E)$,
where $E$ is the exceptional divisor and $k$ is a fixed positive integer.
When $k$ is large the line bundle $\pi^* L^{\otimes k} \otimes {\mathcal O}(-E)$
is ample on $\hat{M}$, so there is an embedding of $\hat{M}$ in a
projective space such that a positive tensor power of this line bundle is
isomorphic to the restriction of the hyperplane line bundle on the
projective space. 
It is proved in Section 3 of \cite{Ki2} that
when $k$ is large enough this linear action satisfies
the following properties:

(i) If $y$ is semistable in $\hat{M}$ then $\pi(y)$ is semistable in $M$.

(ii) If $\pi(y)$ is stable in $M$ then $y$ is stable in $\hat{M}$.

\noindent The rough idea of the proof is to use Proposition 2.1 in conjunction with
the facts that if $k>0$ then the stability and semistability 
with respect to $L^{\otimes k}$ of a point
of $M$  is independent of $k$, and that when $k$ is large the weights
of the action on $H^0(\hat{M},\pi^* L^{\otimes k} \otimes {\mathcal O}(-E))$
of a maximal torus of $G$ can be thought of as small perturbations
of the weights of its action on $H^0(\hat{M},\pi^* L^{\otimes k})$. A similar
argument shows that if $k$ is sufficiently large then the sets $\hat{M}^s$
and $\hat{M}^{ss}$ of stable and semistable points of $\hat{M}$ with respect to this
linearization are independent of $k$.

\begin{rem} \label{newrem41}
Note that the induced symplectic form $\hat{\omega}$ on $\hat{M}$ and 
moment map $\hat{\mu}: \hat{M} \to \lieks$ are of the form
$$\hat{\omega} = k \pi^*\omega + \Omega \mbox{  and  } \hat{\mu} = k \mu \circ \pi + \nu$$
where $\Omega$ and $\nu$ are independent of $k$. Thus if $k \gg 0$ then after
scaling by $1/k$ (which does not change the quotient) the symplectic form and
moment map for $\hat{M}$ are small perturbations $\pi^*\omega + (1/k)\Omega$
and $ \mu \circ \pi + (1/k) \nu$ of the pullbacks to $\hat{M}$ of the symplectic
form and moment map for $M$.
\end{rem}

Now $M$ has semistable points which are not stable if and only if there exists
a nontrivial connected reductive subgroup of $G$ which fixes some
semistable point. If so, let $r>0$ be the maximal dimension of the reductive
subgroups of $G$ fixing semistable points of $M$, and let ${\mathcal R}(r)$
be a set of representatives of conjugacy classes in $G$ of all connected 
reductive subgroups $R$ of dimension $r$ such that
$$Z_R^{ss} = \{x\in M^{ss} : \mbox{$R$ fixes $x$} \}$$
is nonempty. Then 
$$\bigcup_{R \in {\mathcal R}(r)} GZ_R^{ss}$$
is a disjoint union of nonsingular closed subvarieties of $M^{ss}$,
and
$$G Z_R^{ss} \cong G \times_{N^R} Z_R^{ss}$$
where $N^R$ is the normalizer of $R$ in $G$.

By Hironaka's theorem \cite{Hir} we can resolve the singularities of the
closure of $\bigcup_{R \in {\mathcal R}(r)} GZ_R^{ss}$ in $M$ by 
performing a sequence of
blow-ups along nonsingular $G$-invariant closed subvarieties of
$M-M^{ss}$. We then blow up along the proper transform of the
closure of $\bigcup_{R \in {\mathcal R}(r)} GZ_R^{ss}$ to get a nonsingular
projective variety $\hat{M}_1$. The linear action of $G$ on $M$ lifts
to an action on this blow-up $\hat{M}_1$ which can be linearized using
suitable ample line bundles as above,
and it is shown in \cite{Ki2} that the set $\hat{M}_1^{ss}$ of semistable
points of $\hat{M}_1$ with respect to any of these suitable linearizations
of the lifted action is the complement in the inverse image of $M^{ss}$
of the proper transform of the subset
$$\phi^{-1}\left(\phi\left(\bigcup_{R \in {\mathcal R}(r)} 
GZ_R^{ss}\right)\right)$$
of $M^{ss}$, where $\phi:M^{ss} \to \xg$ is the canonical map. Moreover
no point of $\hat{M}_1^{ss}$ is fixed by a reductive subgroup of $G$ of dimension
at least $r$, and a point in $\hat{M}_1^{ss}$ is fixed by a reductive subgroup
$R$ of $G$ of dimension less than $r$ if and only if it belongs to the
proper transform of the subvariety $Z_R^{ss}$ of $M^{ss}$.

The same procedure is now applied in \cite{Ki2} to $\hat{M}_1$ to obtain
$\hat{M}_2$ such that no reductive subgroup of $G$ of dimension at least $r-1$
fixes a point of $\hat{M}_2^{ss}$. After repeating enough times we obtain 
$\tilde{M}$ satisfying $\tilde{M}^{ss}=\tilde{M}^s$, and then the induced map
$\txg \to \xg$ is a partial resolution of singularities.

\begin{rem} \label{oldrem41} If we are only interested in $\tilde{M}^{ss}$ and the partial
resolution $\txg$ of $\xg$, rather than in $\tilde{M}$ itself, then there
is no need in this procedure to resolve the singularities of the closure
of $\bigcup_{R \in {\mathcal R}(r)} GZ_R^{ss}$ in $M$. Instead we can simply blow
$M^{ss}$ up along $\bigcup_{R \in {\mathcal R}(r)} GZ_R^{ss}$ (or equivalently
along each $GZ_R^{ss}$ in turn) and let $\hat{M}_1^{ss}$ be the set of semistable
points in the result.
\end{rem}

\begin{example} \label{egs4} Let $G=SL(2)$ act on $M=(\PP_1)^n$ and suppose
that $n$ is even, so that semistability and stability do not coincide (see
Example 2.4). The semistable elements of $M$ which are fixed by nontrivial
connected reductive subgroups of $G$ are those of the form $(x_1,...,x_n)$
such that there exist distinct $p$ and $q$ in $\PP_1$  
with exactly half of the points $x_1,...,x_n$ equal to $p$ and the rest equal
to $q$. They form 
$$\frac{n!}{2((n/2)!)^2}$$
 $G$-orbits and their stabilizers are all 
conjugate to the maximal torus $T_c = \CC^*$ of $G$. We obtain the 
partial desingularization $\txg$ by blowing up $\xg$ at the points
corresponding to these orbits, or equivalently by blowing up $M^{ss}$
along these orbits, removing the unstable points from the blowup (these
form the proper transform of the set of $(x_1,...,x_n)\in M^{ss}$ such
that exactly half of the points $x_1,...,x_n$ coincide somewhere on $\PP_1$)
and finally quotienting by $G$.
\end{example}

\renorm
\section{Intersection homology and a splitting of the surjection
$\k_M^{ss}:\hk(M^{ss}) \to IH^*(\xg)$}
\label{splitting}

\newcommand{\muchless}{  \ll}
\newcommand{\muchgreater}{ \gg }

In this section, we shall recall the splitting constructed in \cite{Kiem} (see also \cite{KW})
of the surjection $\k_M^{ss}:\hk(M^{ss}) \to IH^*(\xg)$ defined
at (\ref{kss}).

Let $W$ be a (singular) complex projective variety. Then it has
a filtration $W=W_{m}\supseteq
W_{m-1}\supseteq \cdots \supseteq W_0$
by closed subvarieties
which defines a Whitney stratification of $W$ with nonsingular
strata $W_j - W_{j-1}$ of complex dimension $j$, and its intersection cohomology $IH^*(W)$
with complex coefficients and with respect to the middle perversity
can be defined as follows \cite{GM1,GM2}.
Let $IC^{2m-i}(W)$ be the group of chains $\sigma$ of dimension $i$
in $W$ such that
\begin{eqnarray}
dim_{\RR}\,(|\sigma |\cap W_{m-k})\le i-k-1,\\
dim_{\RR}\,(|\partial \sigma|\cap W_{m-k})\le i-k-2.
\end{eqnarray}
Then $IC^*(W)$ is a chain complex whose cohomology is the intersection
cohomology $IH^*(W)$ of $W$. It does not depend
on the choice of the stratification and it is a homeomorphism invariant
\cite{GM1,GM2}. It coincides with ordinary cohomology for
nonsingular varieties, and also for orbifolds
since we are using complex coefficients. If $\sigma \in IC^*(W)$ then neither
$\sigma$ nor $\partial \sigma$ can be contained in $W_{m-1}$, so if $L$ is
any local coefficient system on the nonsingular open subset $W - W_{m-1}$
of $W$ then we can define a chain complex $IC^*(W,L)$ of intersection chains
in $W$ with coefficients in $L$, and thus define the intersection homology
$IH^*(W,L)$ of $W$ with coefficients in $L$.

Any two intersection cohomology classes of complementary degrees in $W$
can be represented by cycles in $W$ which intersect transversely
and only on the
nonsingular part of $W$ at a finite number of points. If
we count these
intersection points with appropriate
signs we obtain a nondegenerate pairing on $IH^*(W)$
which is called the intersection pairing. Thus $IH^*(W)$ satisfies
Poincar\'{e} duality. It also satisfies the properties of the cohomology
of smooth compact K\"ahler manifolds known as the K\"ahler package,
including the existence of a Hodge structure and the hard Lefschetz property.

One of the most useful tools for working with intersection cohomology is the
decomposition theorem of Beilinson, Bernstein, Deligne and Gabber
\cite{BBDG} which tells us that if $f:A\to B$ is a projective map of
complex varieties then there exist closed subvarieties $B_\a$ of $B$
and local systems $L_\a$ on open dense subsets of $B_\a$ such that
\begin{equation} \label{decomp} IH^i(A)\cong
 \bigoplus_{\a} IH^{i-l_{\a}}(B_{\a}, L_{\a})\end{equation}
for suitable integers $l_\a$. If $f$ is birational then there
is some $\a$ such that $B_{\a}=B$ and $L_{\a}=\CC$ and $l_{\a}=0$, so that
$IH^*(B)$ appears as a direct summand of $IH^*(A)$ in this
decomposition.
In particular, the intersection cohomology $IH^*(\xg)$ of
the GIT quotient $\xg$ can be regarded as a direct
summand of the ordinary cohomology $H^*(\xtg)$ of its partial
desingularization $\xtg$, so we get a surjection
\begin{equation} H^*(\xtg) \to IH^*(\xg). \end{equation}
The decomposition (\ref{decomp}) is unfortunately not in
general canonical, but in our situation the hard Lefschetz theorem
can be used to make a canonical choice of decomposition and
hence a canonical choice of surjection
$H^*(\xtg) \to IH^*(\xg)$ (see \cite{Ki3}.)

Our goal is to understand the intersection cohomology of
the singular quotient
$\xg$ in terms of the equivariant cohomology of $M$.
A procedure for computing the intersection
cohomology Betti numbers $\dim IH^j(\xg)$ is given
in \cite{Ki3},  and it was generalized to symplectic
quotients in \cite{Woolf}.
One can compute the equivariant Poincar\'{e} series
$$P_t^K(M^{ss})= \sum_{j \geq 0} t^j \dim H^j_K(M^{ss})$$
of $M^{ss}$ by equivariant Morse theory
applied to the function $|\!|\mu|\!|^2$ as in \cite{Ki1}
(cf. Proposition 2.3 above),
and keep track of
the equivariant Poincar\'{e} series while blowing up until one reaches the
partial desingularization, and then switch to intersection cohomology while
blowing down until one comes back to $\xg$ \cite{KW,Ki3,Woolf}.
The switch is possible since $\tilde{M}^{ss}=\tilde{M}^s$,
so that $\hk(\tilde{M}^{ss})$ is
isomorphic to $H^*(\xtg) = IH^*(\xtg)$.

\begin{example} When $K=SU(2)$ acts on $M=\PP_n$
we have
$$P_t^K(M) = (1+t^2 + t^4 + ... + t^{2n}) (1-t^4)^{-1}$$
and
$$P_t^K(M^{ss}) = P_t^K(M) - \sum_{n/2 < j \leq n}
t^{2(j-1)}(1-t^2)^{-1}$$
(see Example 2.4).
If $n$ is odd so that semistability and stability
coincide then this is a polynomial of degree $2(n-3)$
in $t$ whose coefficients are the (intersection)
Betti numbers of $\xg$. If $n$ is even, then to obtain
the partial desingularization $\xtg$ one must
blow up $M^{ss}$ along the orbit of the element
$[0,...,0,1,0,...,0] \in \PP_n$ corresponding to
the polynomial whose roots in $\PP_1$ are 0 and $\infty$, each with multiplicity $n/2$, and then
remove the unstable points from the blowup. This
gives us
$$P_t(\xtg)= P_t^K(\tilde{M}^{ss}) $$
$$ =
P_t^K(M^{ss}) + (t^2 + t^4 + ... + t^{2(n-3)})(1-t^4)^{-1} -
t^{n-2}(1+t^2+...+t^{n-4})(1-t^2)^{-1}$$
$$=1+2t^2 + 3t^4 + 4t^6 + ... + (\frac{n}{2} -2)t^{n-6}
+(\frac{n}{2}-1)t^{n-4} + $$
$$+ (\frac{n}{2}-1)t^{n-2}
+ (\frac{n}{2} -2)t^{n} + ... + 3t^{2n-10} + 2t^{2n-8} + t^{2n-6}.$$
Finally we study the kernel of the surjection
from $H^*(\xtg)$ to $IH^*(\xg)$ to obtain
the intersection Poincar\'{e} polynomial of $\xg$ as
$$IP_t(\xg) = P_t(\xtg) -
(t^2+t^4 + 2t^6 + 2t^8 + ... +
[(n-2)/{4}]t^{n-4} + $$
$$ + [(n-2)/{4}]t^{n-2} +...+ t^{2n-10} + t^{2n-8}) $$
$$
= 1 + t^2 + 2t^4 + 2t^6 + ... + [n/4]t^{n-4} +
[n/4]t^{n-2} + ... + 2t^{2n-10} + t^{2n-8} + t^{2n-6},$$
where $[a]$ is the integer part of $a$.
For more details see \cite{Ki3}.
\end{example}

The composition of the maps
in the partial desingularization process gives
us a map from $\tilde{M}^{ss}$ to $M^{ss}$ and hence a ring
homomorphism $\hk(M^{ss})\to \hk(\tilde{M}^{ss})$. Via the
decomposition theorem the corresponding maps on quotients induce surjections on intersection cohomology
whose composition
gives us our surjection from $H^*(\xtg) = IH^*(\xtg)$ to $IH^*(\xg)$.
In this way we get
\begin{equation} \k_M^{ss}: \hk(M^{ss})\to
\hk(\tilde{M}^{ss})\cong IH^*(\xtg)\to IH^*(\xg).\end{equation}
This map
$\k_M^{ss}:\hk(M^{ss})\to IH^*(\xg)$ is surjective; the proof
of this in \cite{Ki3} is
flawed but an alternative proof is given in \cite{Woolf2}.

In order to get useful information about the intersection cohomology
$IH^*(\xg)$, a splitting of the map
$\k_M^{ss}: H^*_K(M^{ss}) \to IH^*(\xg)$ was constructed in \cite{Kiem,KW}
under the assumption that the linear action of $G$ on $M$ is weakly balanced in the sense defined below.

\begin{definition}\label{5.3} Suppose a nontrivial connected reductive group
$R$ acts on a vector space $A$ linearly.
Let $\calb$ be the set of the
closest points from the origin to the convex hulls of some weights
of the action. For each $\beta\in \calb$, denote by $n(\b)$ the number
of weights $\a$ such that $\a \cdot\b <\b \cdot\b$. The action
is said to be {\bf weakly linearly balanced} if
 $2n(\b)-2 dim_{\CC} R/BStab\b > dim_{\CC} A-dim_{\CC} R$
for every $\b \in \calb$ where $B$ is a Borel subgroup of $R$.\end{definition}

Let $\mathcal R$ be a set of representatives of the conjugacy classes in $G$
of subgroups which appear as identity components of stabilizers of points
$x\in M^{ss}$ such that $Gx$ is closed in $M^{ss}$. Such subgroups are
always connected reductive subgroups of $G$ (see e.g. \cite{luna}).

\begin{definition}\label{5.4}
Let $G$ be a connected reductive group
acting linearly on a connected nonsingular quasi-projective variety $M$.
The $G$-action is said to be {\bf weakly balanced}
 if for each $R\in \mathcal{R}$ the linear action of $R$ on the normal space
$\caln_x$ at any $x\in Z^s_R$ to $GZ^{ss}_R$ is
weakly linearly balanced and so is the action
of $(R\cap N^{gR'g^{-1}})/gR'g^{-1}$ on the $gR'g^{-1}$-fixed
linear subspace $\mathcal{N}_x^{gR'g^{-1}} = $
$Z^{ss}_{gR'g^{-1}}\cap \caln_x$
for each $R'\in \mathcal{R}$ satisfying $gR'g^{-1}
\subseteq R$. \end{definition}

For example, a $\CC^*$ action
on $\PP^n$ is weakly balanced if and only if the number (counting multiplicities) of
positive weights is same as the number of negative weights. The actions
described in Examples 2.4 and 2.5 of $SL(2)$ and its maximal torus
$\CC^*$ on $\PP_n$ and $(\PP_1)^n$ are weakly balanced. More examples
are provided by the
(compactified) moduli spaces
of holomorphic vector bundles of any rank and degree
over a fixed Riemann surface (see Example 3.5 of
\cite{Kiem}).

For $R\in \mathcal R$,
we consider the natural map
\begin{equation}  G \times_{N^R}Z^{ss}_R\to GZ^{ss}_R
\end{equation}
where $N^R$ is the normalizer of $R$ in $G$, and the corresponding map
\begin{equation}   \hk(GZ^{ss}_R)\to \hk(G\times_{N^R}Z^{ss}_R)\cong
H^*_{N^R}(Z^{ss}_R)
\cong [H^*_{N^R_0/R}(Z^{ss}_R)\otimes H^*_R)]^{\pi_0N^R}\end{equation}
where the subscript $0$ means the identity component.
For any $\zeta\in H^*_K(M^{ss})$ we let $\zeta|_{G\times _{N^R} Z^{ss}_R}$
denote the image of $\zeta$ under the composition of the above map and
the restriction map $H^*_K(M^{ss})\rightarrow H^*_K(GZ^{ss}_R)$.
Then the main result of \cite{Kiem} is the following splitting
of the map
$\k_M^{ss}: \hk(M^{ss}) \to IH^*(\xg)$ obtained by ``truncating along each stratum'' when the action is
weakly balanced.
\begin{theorem}\label{5.5} \cite{Kiem} Let
$$  V_M=\{\zeta\in H^*_K(M^{ss}):\,\,\zeta|_{G\times _{N^R} Z^{ss}_R}
\in [\oplus_{i<n_R}H^*_{N_0^R/R}(Z^{ss}_R)\otimes H^{i}_R]^{\pi_0N^R}
\mbox{ for  each } R\in \mathcal R\}$$
where $n_R=dim_{\CC}\caln_x-dim_{\CC} R=dim_{\CC} \caln_x/\!/R$ and
$\caln_x$ is the normal space to $GZ_R^{ss}$
at any $x\in Z^s_R$. If the action of $G$ is weakly balanced, then
the restriction
$$\k_M^{ss}:V_M\to IH^*(\xg)$$
of the map $\k_M^{ss}:\hk(M^{ss}) \to IH^*(\xg)$
is an isomorphism.\end{theorem}

\begin{example}\label{5.6}
We continue Example \ref{5.2}. In the terminology of
Theorem \ref{5.5} we have $n_R=5$, and so we have to remove
$\CC\{\xi^i\rho^j\,:\, i=0,1,\,j\ge 3\}$ to get $V_M$. Hence,
$$V=\oplus_{0\le i\le 6} V^{2i}$$
where
$$V^0={\CC},\,\,\,\,\, V^2={\CC}\{\rho,\xi\},\,\,\,\,\,\,
V^4={\CC}\{\rho^2,\xi\rho, \xi^2\},$$
$$V^6={\CC}\{\xi\rho^2,\xi^2\rho,\xi^3\},\,\,\,\,\,
V^8={\CC}\{\xi^2\rho^2,\xi^3\rho,\xi^4\},$$
$$V^{10}={\CC}\{\xi^2\rho^3,\xi^3\rho^2\},\,\,\,\,\,\,
V^{12}={\CC}\{\xi^2\rho^4\}.$$
Therefore, the intersection Poincar\'{e} series for $\PP^7/\!/\CC^*$ is
$$1+2t^2+3t^4+3t^6+3t^6+3t^8+2t^{10}+t^{12}.$$
This formula could also be obtained by the sort
of calculation described for $\PP_n /\!/SL(2)$
in Example 5.1, but such a calculation is usually lengthier.
\end{example}

\begin{rem} The equivariant cohomology of
$M$ and of $M^{ss}$ and the intersection cohomology
of $\xg$ carry natural Hodge structures, and the
maps we have been considering (in particular $\k_M$
and $\k_M^{ss}$) respect these Hodge structures. Thus we can use these methods to calculate Hodge numbers
as well as Betti numbers. In particular if every
cohomology class of $M$ is of Hodge type
$(p,p)$ for some $p$
then the same is true for $\xg$, as can be seen in the last example
(cf. \cite{Ki1} Section 14).
\end{rem}

\renorm
\section{Intersection pairings via the splitting $V_M$}
\label{pairings via the splitting}

In this section, we study the intersection
pairing in $IH^*(\xg)$ via $V_M$. Throughout this section, we assume that the
action of $G$ is weakly balanced and thus we have the isomorphism
$$\k_M^{ss}:V_M\to IH^*(\xg).$$

Let $\tau$ be the top degree class in $V_M$ corresponding to the
fundamental class in $IH^m(\xg)$ where $m$ is the real dimension of
the quotient $\xg$. It comes from a class in
the compactly supported cohomology group
$H^m_c(M^s/G)$ which can be
represented by a closed differential form
with compact support, via
the composition
$$H^m_c(M^s/G)\cong H^m(\xg,M^{sss}/\!/G)\cong
H^m_K(M^{ss},M^{sss})\to H_K^m(M^{ss}).$$
Then we have the following theorem \cite{Kiem}.
\begin{theorem}\label{6.1} Let $\a, \,\b$ be two classes of complementary
degrees in $V_M$ with respect to $m$. Then their product $\a \b$ in $\hk(M^{ss})$
is the top degree class $\tau$ in $V_M$ multiplied by the scalar
\begin{equation}  \langle \k_M^{ss}(\a),\k_M^{ss}(\b) \rangle_{IH^*(\xg)}
\end{equation}
where $\langle \cdot,\cdot\rangle _{IH^*(\xg)}$ denotes 
 the intersection pairing
in $IH^*(\xg)$.\end{theorem}

\begin{example}\label{6.2} We continue
Examples 2.5 and \ref{5.6}. Notice that by the
Gr\"obner basis in Example 2.5,
$\xi^3\rho^3=0$, $\xi^4\rho^2=-\frac13
\xi^2\rho^4$, $\xi^6=\xi^2\rho^4$, and $\xi^5\rho=0$. Hence, for example,
the matrix for the pairing
$V^6\otimes V^6\to \CC$
is up to a constant
$$\left( \begin{array}{ccc}
1 &0 &-\frac13\\
0 &-\frac13 &0\\
-\frac13 &0 &1
\end{array} \right)$$
One can similarly compute the pairings for other classes.
\end{example}

Since the map $\hk(M^{ss})\to \hk(\tilde{M}^{ss})\cong H^*(\tilde{M}/\!/G)$
induced from the maps in the partial desingularization process is a ring
homomorphism, the relation $$\a \b =\langle \k_M^{ss}(\a),\k_M^{ss}(\b)\rangle _{IH^*(\xg)}
\,\,\,\tau$$
in Theorem 6.1 is preserved. The intersection pairing
$\langle \k_M^{ss}(\a),\k_M^{ss}(\b)\rangle _{IH^*(\xg)}$ therefore
equals the evaluation $\k_{\tilde{M}}^{ss}(\a \b)[\tilde{M}/\!/G]$ of
the image in $H^*(\xtg)$ of the product $\a \b\in \hk(M^{ss})$
on the fundamental class $[\tilde{M}/\!/G]$,
because $\k_{\tilde{M}}^{ss}(\tau)[\tilde{M}/\!/G]=1$.
Hence, we get the following.

\begin{proposition}\label{6.3} Let $\a, \b$ be classes of complementary degrees
in $V_M$. Then
$$ \langle \k_M^{ss}(\a),\k_M^{ss}(\b)\rangle _{IH^*(\xg)}=\k_{\tilde{M}}^{ss}(\a \b) [\tilde{M}/\!/G].$$
\end{proposition}

\begin{rem} Via the natural map $H^*(\xg) \to IH^*(\xg)$
the pairing in $IH^*(\xg)$ determines the intersection pairing in
the ordinary cohomology $H^*(\xg)$, which may be degenerate for singular
quotients. The quotient map $EK\times_K M^{ss}\to \xg$
induces a ring homomorphism $H^*(\xg)\to \hk(M^{ss})$ which factors
through $V_M$, and the composition $H^*(\xg)\to \hk(M^{ss})\to IH^*(\xg)$
is the natural map $H^*(\xg) \to IH^*(\xg)$
which preserves the pairing.
\end{rem}

\renorm
\section{Pairings in intersection cohomology}

In this section and the next we consider pairings in $IH^*(\xg)$ and pairings in the cohomology
$H^*(\tilde M \git G)$ of classes in the image of the composition
$$H^*_K(M) \rightarrow H^*_K(\tilde M) \rightarrow H^*(\tilde M \git G)$$
of the pullback from $M$ to $\tilde{M}$ and the map $\kappa_{\tilde M}$.
As in previous sections we abuse notation by suppressing the pullback and
writing $\kappa_{\tilde M}(\alpha)$ for the image
of $\alpha\in H^*_K(M)$. Since these maps are ring homomorphisms,
such a pairing is given simply by evaluating the product against the
fundamental class; i.e.
$$\langle \kappa_{\tilde M}(\alpha) , \kappa_{\tilde M}(\beta)
\rangle_{H^*(\tilde M \git G)} = \kappa_{\tilde M}(\alpha\beta)
[\tilde M \git G].$$
Furthermore in  Section   \ref{pairings via the splitting} it was
shown that if the $G$ action is weakly balanced and
$\alpha|_{M^{ss}}$ and $\beta|_{M^{ss}}$ lie in the subspace
$V_M$ of $\hk(M^{ss})$ which is isomorphic to $IH^*(\xg)$ then
$$\langle \kappa_{\tilde M}(\alpha) , \kappa_{\tilde M}(\beta)
\rangle_{IH^*( M \git G)} = \kappa_{\tilde M}(\alpha\beta)[\tilde M
\git G].$$
In this section we will find a formula for
$\k_{\tilde{M}}(\alpha\beta)[\txg]$
in this special case; in
the next section we will study evaluations of the form
$\k_{\tilde{M}}(\alpha\beta)[\txg]$
for any $\alpha, \beta \in H^*_K(M)$.

Recall that we are assuming that $M$ is a nonsingular complex
projective variety embedded in a complex projective space
$\PP_n$, and that $G$ acts on $M$ via a complex representation
$\rho:G \to GL(n+1)$ of $G$ such that $\rho(K) \subseteq  U(n+1)$.
This representation $\rho$ gives us a lift of the action of $G$ on
$M$ to the
hyperplane line bundle over $M$, i.e. a linearization of the
action of $G$ on $M$. We can change the linearization without
changing the action of $G$ on $M$ by multiplying $\rho$ by any
character
$\chi:G \to \CC^*$ of $G$. If we
identify $\chi$ with an element of $\lieks$ in the usual way
by taking the derivative at the identity of the restriction of
$\chi$ to $K$, then this change in linearization corresponds
to shifting the moment map
$\mu:M \to \lieks$ by the central element $\chi$ of $\lieks$.
The quotient
$$M/\!/_{\chi}G$$
of $M$ by $G$ with respect to this shifted linearization can
be identified topologically with the quotient of the Zariski open subset
$$M_{\chi}^{ss} = \{m \in M | \chi \in \mu(\overline{G \cdot m}) \}$$
of $M$ by the equivalence relation $\sim$ such that $x \sim y$
if and only if
$$\overline{G \cdot x} \cap \overline{G \cdot y} \cap M_{\chi}^{ss}
\neq \emptyset.$$
Just as we have a homeomorphism $M \git G \cong \mu^{-1}(0)/K$,
so we have homeomorphisms
$$M \git_{\chi} G \cong \mu^{-1}(\chi)/K$$
for any such $\chi$. Moreover we can generalize the construction
of $M_{\chi}^{ss}$ and $M \git_{\chi} G \cong \mu^{-1}(\chi)/K$
to any central element $\chi $ of $\lieks$. In particular for
any $\xi \in \liets$ we can define
\begin{equation} \label{new3} M_{\xi,T}^{ss} = \{m \in M | \xi \in \mu_T(\overline{T_c \cdot m})\} \end{equation}
and its quotient $M \git_{\xi} T_c = M_{\xi,T}^{ss}/ \sim$ which is homeomorphic to
$ \mu_T^{-1}(\xi)/T$. When $\xi$ is a regular value of $\mu_T$ then
just as at (1.1) we get a map
$$H^*_{T_c}(M) \rightarrow H^*_{T_c}(M_{\xi,T}^{ss})
\cong H^*(M \git_\xi T_c)$$
which we shall denote by $\kappa^T_{M,\xi}$.

By the convexity theorem of Atiyah \cite{Aconv} and Guillemin and Sternberg \cite{GS}
the image $\mu_T(M)$ of the moment map $\mu_T$ is a convex polytope
in $\liet^*$ and the dense set of regular values
 is the disjoint union of finitely
many open convex subpolytopes
$$\Delta_1 \cup \ldots \cup \Delta_r.$$ In fact $\mu_T(M)$ is
the convex
hull in $\liets$ of the finite set $\{\mu_T(F):F \in \calf\}$, and it is divided
by walls of codimension one into subpolytopes which are convex
hulls of subsets of  $\{\mu_T(F):F \in \calf\}$ and whose interiors consist
of regular values of $\mu_T$.
Here $\calf$ is the set of components of the fixed point set of the
action of $T$ on $M$.

Suppose $\xi \in \Delta_i$ and $\zeta \in \overline{\Delta_i}$.
Then $M_\xi^{ss} \subseteq M_\zeta^{ss}$ and this inclusion
induces a birational map
$$\mu_T^{-1}(\xi)/T \cong M_\xi^{ss} \git T_c \rightarrow
M_\zeta^{ss} \git T_c\cong\mu_T^{-1}(\zeta)/T.$$
Of course if $\zeta \in \Delta_i$ then this map is an isomorphism\footnote{Each
face of $\Delta_i$ is similarly divided into subpolytopes, and if $\xi$ and $\zeta$
both lie in the interior of the same subpolytope in a face of $\Delta_i$ then
$$\mu_T^{-1}(\xi)/T \cong M_\xi^{ss} \git T_c =
M_\zeta^{ss} \git T_c\cong\mu_T^{-1}(\zeta)/T.$$},
but this is not true in general when $\zeta$ lies in the boundary of $\Delta_i$.

\begin{proposition}
Suppose the action of $G$ on $M$ is weakly balanced.
For $\alpha,\beta \in H^*_K(M)$ such that
$\alpha|_{M^{ss}} ,\beta|_{M^{ss}} \in V_M$ we have
$$\kappa_{\tilde M}(\alpha \beta)[\tilde M \git G] =
\frac{n_0 (-1)^{n_+}}{n_0^T |W|} \kappa^T_{M,\xi}(\alpha\beta
\cald^2)[M \git_\xi T_c]$$
for any $\xi \in \Delta_i$ with $0 \in \overline{\Delta_i}$.
\end{proposition}

\noindent{\em Proof:} Recall from Theorem \ref{6.1} and Proposition \ref{6.3} that  if $\alpha$ and $\beta$ are of complementary degree then
\begin{equation}\label{tau equation}\alpha\beta|_{M^{ss}}
= \kappa_{\tilde M}(\alpha \beta)[\tilde M \git G] \tau
\end{equation}
where $\tau \in H^{\dim M \git G}_K(M^{ss})$ is the class in $V_M$
corresponding to the dual in $IH^{\dim M \git G}(M \git G)$ of the
class of a point in $IH_0(M \git G)$.

Let $U$ be open in the subset $M^s \subseteq M^{ss}$ of stable
points in $M$
for the $G$ action (and recall that $U$ is then contained in the
stable part of $M$ with respect to the $T_c$ action with linearization
induced from that of $G$). Applying to compactly supported cohomology on $U$
the arguments used by Martin in
\cite{Martin,Martin2} (see also \cite{arbrank} Section 3) to prove
(\ref{ktot}),  we find
that $\tau \cald^2 \in H^*_{T_c}(M_{0,T}^{ss})$ is the image
under the composition$$H^{\dim M \git T_c}_c(U / T_c) \cong
H^{\dim M \git T_c}_{T_c,c}(U)
\rightarrow H^{\dim M \git T_c}_{T_c}(M_{0,T}^{ss})$$
of the unique generator $\gamma$ of $H^{\dim M \git T_c}_c(U / T_c)$
with the normalization
$$\gamma [ U / T_c] = \frac{n_0^T |W|}{n_0 (-1)^{n_+}} .$$
Let us now take $U = M_{\xi,T}^{ss} \cap M^s$. There is then a
commutative diagram
\[\xymatrix{& U / T_c \ar[dl] \ar[dr] & \\M \git_\xi T_c \ar[rr] && M \git T_c.}\]
By equation (\ref{tau equation}) we have
$$\alpha\beta\cald^2|_{M^{ss}} = \kappa_{\tilde M}(\alpha \beta)
[\tilde M \git G] \tau\cald^2|_{M^{ss}}.$$
So, using the above diagram and description of $\tau \cald^2$,
we see that
\begin{eqnarray*}
                    \frac{n_0 (-1)^{n_+}}{n_0^T |W|} \kappa_{M,p}^T(\alpha\beta\cald^2)
[M \git_\xi T_c] &=&\kappa_{\tilde M}(\alpha\beta)[\tilde M \git G] \left(
 \frac{n_0 (-1)^{n_+} \gamma}{n_0^T |W|} 
[U/T_c] \right) \\&=
& \kappa_{\tilde M}(\alpha\beta)[\tilde M \git G]
\end{eqnarray*}
as required.

\begin{rem}This result can be viewed as an extension of
the well known fact that for $\xi,\zeta \in \Delta_i$ and any $\eta \in
H^*_{T_c}(M)$ we have
$$\kappa_{M,\xi}^T(\eta)[M \git_\xi T_c] = \kappa_{M,\zeta}^T(\eta)[M
\git_\zeta T_c].$$
Provided the action is weakly balanced with respect to $\zeta$\footnote{By ``weakly
balanced with respect to $\zeta$'' we mean weakly balanced with respect to the
moment map $\mu - \zeta$, in other words we subtract
$\zeta$ from all the weights.} and
we restrict $\eta$ to be a scalar multiple of $\tau \cald^2$ then the same
formula is valid when $\zeta \in \overline{\Delta_i} \setminus \Delta_i$.
\end{rem}

Since $\xi$ is a regular value of $\mu_T$ we can evaluate
$$\kappa_{M,\xi}^T(\alpha\beta\cald^2)[M \git_\xi T_c]$$
using a residue formula as at (\ref{res}) which involves a sum over the set
$\mathcal{F}$ of fixed point components of the $T$ action on $M$.
Together with the results of  Section  \ref{pairings via the splitting}
this allows us, when the $G$ action
is weakly balanced, to write down residue formulae for computing
pairings in
the intersection cohomology $IH^*(M \git G)$, as follows.

\begin{theorem} Suppose that the $G$ action on $M$ is
weakly balanced and that $\a, \b \in \hk(M)$ have degrees
whose sum is the real dimension of $\xg$. Suppose also that
$\alpha|_{M^{ss}}$  and
$\beta|_{M^{ss}}$ lie in the subspace $V_M$ of $\hk(M^{ss})$
which is isomorphic to $IH^*(\xg)$. Then
$$\langle \kappa_M(\alpha) , \kappa_M(\beta) \rangle_{IH^*( M \git G)}
 = \kappa_{\tilde M}(\alpha\beta)[\tilde M \git G]$$$$= \frac{n_0
(-1)^{s+n_+}}{|W| \vol (T)} \res \bigl( {\mathcal{D}}^2\sum_{{F} \in
{\calf}} \int_{F}\frac{i_{F}^*(\a\b e^{{(\bom - \e)}})(\xvec)}{e_{{F}}}
[d\xvec]\bigr)$$
for any sufficiently small $\e \in \liets$ which is a regular value of
the moment map $\mu_T$.
\end{theorem}

\begin{example} \label{caps}
Let us consider again the action of $G=SL(2)$ on $M=\PP_n$
(Example 2.4 and 5.1).
 The equivariant
cohomology $H^*_T(M)$ of $M$ with respect to the maximal torus $T$ of
$K=SU(2)$ is generated by two equivariant cohomology classes $\xi$ and $\zeta$,
both of degree two, subject to the relation
$$(\xi - n\zeta)(\xi - (n-2)\zeta)...(\xi +(n-2)\zeta)(\xi +n\zeta)=0,$$
where the nonidentity element of the Weyl group $W$ of $K$ sends
$\zeta$ to $-\zeta$ and fixes $\xi$. Thus the equivariant cohomology
$H^*_K(M)$ of $M$ with respect to $K$ is generated
to $\xi$ and $\zeta^2$ subject to the same relation. The
fixed point sets for the action of $T$ are the $n+1$ points
represented by the weight vectors; the values taken by the moment
map on the fixed points are just the weights
$-n,2-n,...,n-2,n$ up to a universal scalar multiple. If $n$ is odd
then semistability equals stability, and if $\eta\in H^*_K(M)$
has degree $2(n-3)=\dim_{\RR}(\xg)$ and is given by a polynomial
$q(\xi,\zeta^2)$ in the generators $\xi$ and $\zeta^2$, then
the residue formulas
(1.4) or (3.7) give us
$$\k_M(\eta)[\xg] = \res_{X=0}\bigl( 4X^2 \sum_{j:0<n-2j\leq n}
\frac{q((n-2j)X,X^2)}{\prod_{k\neq j}((n-2k)X - (n-2j)X)}\bigr)$$
$$= \res_{X=0}\bigl( \sum_{j:0<n-2j\leq n}
\frac{q((n-2j)X,X^2)}{(2X)^{n-2} \prod_{k\neq j}(k-j)}\bigr)$$
(cf. Example 3.2). If, on the other hand,
$n$ is even then  Theorem
7.3 shows us that if  $\a, \b \in \hk(M)$ satisfy $\deg(\a) + \deg(\b) =
2(n-3)$ and if $\alpha|_{M^{ss}}$ and $\beta|_{M^{ss}}$ lie in the
subspace $V_M$ of $\hk(M^{ss})$ then
$$\langle \kappa_M(\alpha) , \kappa_M(\beta) \rangle_{IH^*( M \git G)}
 = \kappa_{\tilde M}(\alpha\beta)[\tilde M \git G]$$
$$= \res_{X=0} \sum_{j: 0<n-2j\leq n} \frac{q((n-2j)X,X^2)}{(2X)^{n-2}
\prod_{k \neq j} (k-j)}$$
where $\a\b = q(\xi, \zeta^2)$.
\end{example}

\renorm
\section{Pairings on the partial desingularization}

\newcommand{\quott}{/\! /}

In this section we consider pairings $ \kappa_{{\tilde{M}}}(\a\b)[{\tilde{M}}/\!/G]$ in the cohomology
$H^*(\tilde M \git G)$ of classes $ \kappa_{{\tilde{M}}}(\a)$ and
$ \kappa_{{\tilde{M}}}(\b)$  in the image of the composition
$$H^*_K(M) \rightarrow H^*_K(\tilde M) \rightarrow H^*(\tilde M \git G)$$
of the pullback from $M$ to $\tilde{M}$ and the map $\kappa_{\tilde M}$.

First note that applying (\ref{ktot}) to $\tilde{M}$ we have
\begin{equation} \label{new7} \kappa_{{\tilde{M}}}(\a\b)[{\tilde{M}}/\!/G] = \frac{n_0
(-1)^{n_+} }{n_0^T |W|} \kappa_{{\tilde{M}}}^T(\a \b
{\mathcal{D}}^2)[{\tilde{M}}/\!/T_c]. \end{equation}

\begin{rem} \label{new1} As at (\ref{ktot}) (cf. Theorem 3.1(ii)) we can replace the pairing 
$\kappa_{{\tilde{M}}}^T(\a \b
{\mathcal{D}}^2)[{\tilde{M}}/\!/T_c]$ by the evaluation on $[\tilde{\mu}^{-1}(0)/T]$
of the cohomology class induced by $(-1)^{n_+} \a\b\cald$ if we want to
make sure that we are only working with the semistable part $\tilde{M}^{ss}$
of $\tilde{M}$.
Indeed we can think of $(-1)^{n_+} \cald$ as representing the Poincar\'{e} dual
to $\tilde{\mu}^{-1}(0)/T$ in $\tilde{\mu}_T^{-1}(0)/T = \tilde{M}/\!/T_c$ or the equivariant
Poincar\'{e} dual to $\tilde{\mu}^{-1}(0)$ in $\tilde{\mu}_T^{-1}(0)$.
\end{rem}

Now let $\xi$ be any regular value of the $T$-moment map $\mu_T$ for $M$. Let $\tilde{M}^{ss}_{\xi,T}$
be defined as at (\ref{new3}); then it follows from  Section  3 of \cite{Ki2} that if we choose $k \gg 0$
and if $x \in \tilde{M}$ lies in $\tilde{M}^{ss}_{\xi,T}$
then its image $\pi (x)$ in $M$ lies in $M^{ss}_{\xi,T}$. Thus $\pi:\tilde{M} \to M$ induces a birational
morphism
$$\pi_\xi : \tilde{M}/\!/_\xi T_c \to M/\!/_\xi T_c$$
and we have $(\pi_\xi)_*[\tilde{M}/\!/_\xi T_c] = [M/\!/_\xi T_c]$. Since
$\pi_\xi : H^*({M}/\!/_\xi T_c) \to H^*(\tilde{M}/\!/_\xi T_c)$ fits into the commutative diagram
\begin{equation} \label{diag7.2}
\xymatrix{
& {H^*_{T}(M)} \ar[dr] \ar[dl] & \\
{H^*_{T}(\tilde{M})} \ar[dr] & &
{H^*(M/\!/_{\xi}T_c )} \ar[dl] \\
& {H^*({\tilde{M}}/\!/_{\xi}T_c )}
}
\end{equation}
this implies that
\begin{equation} \label{new2}
\kappa_{{\tilde{M}},\xi}^T(\a \b
{\mathcal{D}}^2)[{\tilde{M}}/\!/_\xi T_c] = \kappa_{{{M},\xi}}^T(\a \b
{\mathcal{D}}^2)[{{M}}/\!/_\xi T_c],
\end{equation}
and since $\xi$ is a regular value of $\mu_T$, the residue formula Theorem 3.1 gives us

\begin{lemma} \label{new6}  If the degrees of $\a, \b \in \hk(M)$ add
up to the real dimension of $\xg$, then
$$
\kappa^T_{{{M,\xi}}}(\a \b \cald^2)[{{M}}/\!/_\xi T_c] = \frac{n^T_0
(-1)^{\ell}}{ \vol (T)} \res \bigl( {\mathcal{D}}^2
\sum_{{F} \in {\calf}} \int_{{F}}
\frac{i_{{F}}^*(\a\b e^{{\bom}-\xi})(\xvec)}{e_{{F}}}
[d\xvec]\bigr).$$
\end{lemma}

This means that to calculate the pairing $
\kappa_{{\tilde{M}}}(\a \b)[{\tilde{M}}/\!/G] $
it suffices to calculate the difference between
\begin{equation} \label{anew}
\kappa_{{\tilde{M}},\xi}^T(\a \b
{\mathcal{D}}^2)[{\tilde{M}}/\!/_\xi T_c]
\end{equation}
and
\begin{equation} \label{bnew}
\kappa_{{\tilde{M}}}^T(\a \b
{\mathcal{D}}^2)[{\tilde{M}}/\!/ T_c] .
\end{equation}
(Note that (\ref{bnew}) is just the special case of (\ref{anew})
when $\xi=0$.) This difference can be calculated using the version
of nonabelian localization due to Guillemin and Kalkman \cite{GK} and
independently to Martin \cite{Martin, Martin2}.

In fact since the construction of $\tilde{M}^{ss}$ from $M^{ss}$ described in
Section 4 takes place in stages, it is actually easier to consider a
single stage of the construction, when $\hat{M}$ is obtained by blowing $M$ up along the closure
of $G Z_R^{ss}$ for a suitable reductive subgroup $R$ of $G$, after first
resolving the singularities of $\overline{GZ_R^{ss}}$ (or equivalently, if we are
only interested in the semistable points of $\hat{M}$, we can simply blow $M^{ss}$ up
along $G Z_R^{ss}$: see Remark \ref{oldrem41}).
The argument which gave (\ref{new2}) also gives us
\begin{equation} \label{new4}
\kappa_{{\hat{M}},\xi}^T(\a \b
{\mathcal{D}}^2)[{\hat{M}}/\!/_\xi T_c] = \kappa_{{{M}}}^T(\a \b
{\mathcal{D}}^2)[{{M}}/\!/_{\x}T_c]
\end{equation}
when $\xi$ is a regular value of $\mu_T$ and $k$ is chosen sufficiently large (depending
on $\xi$). Let us choose $\xi \in \liets$  to lie in a connected component
$\Delta_i$ of the set of regular values of $\mu_T$ for which $0 \in \bar{\Delta}_i$.
Let $\hat{\mu}$ and $\hat{\mu}_T$ be the moment maps for the actions of $K$ and
$T$ on $\hat{M}$, and choose $\hat{\xi} \in \liets$ to lie in the intersection of $\Delta_i$
and a connected component of the
set of regular values of $\hat{\mu}_T$ which contains $0$ in its closure. Because the
choice of $k$, and hence also of the moment maps $\hat{\mu}$ and $\hat{\mu}_T$, depends on
$\xi$, we cannot necessarily choose $\hat{\xi} = \xi$. However it suffices to
calculate the difference
\begin{equation} \label{new5}
\kappa_{{\hat{M}},\xi}^T(\a \b
{\mathcal{D}}^2)[{\hat{M}}/\!/_\xi T_c] - \kappa_{{\hat{M},\hat{\xi}}}^T(\a \b
{\mathcal{D}}^2)[{\hat{M}}/\!/_{\hat{\xi}}T_c],
\end{equation}
since combining this with (\ref{new2}) and
(\ref{new4}) and iterating the calculation will give us the
difference between $\kappa_{{{M},\xi}}^T(\a \b
{\mathcal{D}}^2)[{{M}}/\!/_{\xi} T_c]$
and
$\kappa_{{\tilde{M}},\tilde{\xi}}^T(\a \b
{\mathcal{D}}^2)[{\tilde{M}}/\!/_{\tilde{\xi}} T_c]$
for any $\tilde{\xi}$ in a connected component of the set of regular values of
$\tilde{\mu}_T$ which contains 0 in its closure. As 0 is itself a regular value
of $\tilde{\mu}_T$, we can choose $\tilde{\xi}$ to be 0 and use Lemma
\ref{new6} to calculate the pairing (\ref{new7}) which is our goal.

\begin{rem} \label{new9}
Notice also that since $\xi \in \Delta_i$ with $0 \in \bar{\Delta}_i$, we have
$$M^{ss}_{\xi,T} \subseteq  M^{ss}_{0,T}$$
where $M^{ss}_{0,T}$ retracts $T$-equivariantly onto $\mu_T^{-1}(0)$, and
similarly
$$\hat{M}^{ss}_{\hat{\xi},T} \subseteq  \hat{M}^{ss}_{0,T}$$
where  $\hat{M}^{ss}_{0,T}$ retracts $T$-equivariantly onto $\hat{\mu}_T^{-1}(0)$.
Thus by Remark \ref{new1} we do not need to resolve the singularities of the closure
of $GZ_R^{ss}$ and construct the whole of $\hat{M}$ and $\tilde{M}$ in order to
carry out these constructions; it suffices to consider the blow-up of $M^{ss}$ along
$GZ_R^{ss}$ and iterations of this process.
\end{rem}

We know that the image $\hat{\mu}_T(\hat{M})$ of the moment map $\hat{\mu}_T$ is a
convex polytope which is divided by walls of codimension one into subpolytopes
whose interiors consist of regular values of $\hat{\mu}_T$.
If $\eta \in \hT (\hat{M})$ there
is no change in
$\kappa_{{\hat{M}},\xi}^T(\eta )[{\hat{M}}/\!/_\xi T_c] $ as $\xi$ varies within a connected
component of the set of regular values of $\hat{\mu}_T$, so it is enough to be able to
calculate the change in
$\kappa_{{\hat{M}},\xi}^T(\eta )[{\hat{M}}/\!/_\xi T_c] $ as $\xi$ crosses a wall of codimension
one. Any such wall is of the form
$$\hat{\mu}_T ( \hat{M}_1)$$
where $\hat{M}_1$ is a connected component of the fixed point set in $\hat{M}$ of a circle
subgroup $T_1$ of $T$. The quotient group $T/T_1$ acts on $\hat{M}_1$, which is a symplectic
submanifold of $\hat{M}$, and the restriction of the moment map $\hat{\mu}_T$ to $\hat{M}_1$
has an orthogonal decomposition
$$\hat{\mu}_{T/T_1} \oplus \hat{\mu}_{T_1}$$
where $\hat{\mu}_{T/T_1}: \hat{M}_1 \to (\liet / \liet_1)^*$ is a moment map for the
action of $T/T_1$ on $\hat{M}_1$, and $\mu_{T_1}$ is constant because the action of $T_1$ on
$\hat{M}_1$ is trivial. Guillemin and Kalkman \cite{GK} show that the change in
$\kappa_{{\hat{M}},\xi}^T(\eta )[{\hat{M}}/\!/_\xi T_c] $ as $\xi$ crosses a section of the wall
$\hat{\mu}_T(\hat{M}_1)$ whose orthogonal projection onto $(\liet /\liet_1)^*$ contains
a regular value $\xi_1$ for $\hat{\mu}_{T/T_1}$ is given by
\begin{equation} \label{new8}
\kappa^{T/T_1}_{\hat{M}_1,\xi_1}( \res _{\hat{M}_1}(\eta)) [ \hat{M}_1 /\!/_{\xi_1} (T/T_1)_c].
\end{equation}
Here the residue operation
$$\res _{\hat{M}_1} : \hT ( \hat{M}) \to H^*_{T/T_1} (\hat{M}_1 )$$
is obtained by choosing a coordinate system $X = (X_1,...,X_{\ell})$ for $\liets$ such that
$X_1$ defines an integer basis for $\liet_1^*$ and $(X_2,...,X_{\ell})$ is a coordinate
system for the dual of a Lie algebra of a codimension one subtorus of $T$ whose
intersection with $T_1$ is finite, and setting
$$\res _{\hat{M}_1} (\eta) = \res_{X_1 =0} (\frac{\eta|_{\hat{M}_1}}{e_{\hat{M}_1}})$$
where $e_{\hat{M}_1}$ is the $T$-equivariant Euler class of the normal bundle to
$\hat{M}_1$ in $\hat{M}$. We make sense of this residue as an element of
$H^*_{T/T_1} (\hat{M}_1 )$ by using $(X_2,...,X_{\ell})$ as coordinates on $(\liet/\liet_1)^*$
and writing
$$ e_{\hat{M}_1} = \prod_{i=1}^{\mbox{codim}\hat{M}_1} (m_i X_1 + e_i)$$
where each $m_i \in \ZZ \setminus \{0\}$  is a weight for the action of $T_1$ on the
normal bundle and each $e_i$ can be identified with an element of
$H^*_{T/T_1} (\hat{M}_1 )$.

\begin{lemma} \label{l8.4}
Let $\xi \in \liets$ lie in a connected component $\Delta_i$ of the set of
regular values of $\mu_T$ for which $0 \in \bar{\Delta_i}$, and
let $\hat{\xi} \in \liets$ lie in the intersection
of $\Delta_i$ and a connected component of the set of regular values
of $\hat{\mu}_T$ which contains $0$ in its centre.
Then in order to calculate the difference (\ref{new5})
the only wall crossing terms we need to consider correspond
to components
$\hat{M_1}$ of fixed point sets of circle subgroups $T_1$ satisfying
$$\emptyset \ne \pi(\hat{M_1}) \cap M^{ss} \subseteq G Z_R^{ss}$$
so that $\hat{M_1}$ is contained in the exceptional divisor of $\pi:
\hat{M} \to M$.
\end{lemma}

\Proof ~From Remark \ref{newrem41} we see that by choosing $k$
sufficiently large we can
assume that (after scaling by $1/k$) the wall
$\hat{\mu}_T(\hat{M}_1)$ is contained in an arbitrarily small
neighbourhood of $\mu_T( \pi(\hat{M}_1))$. Notice that
$\pi (\hat{M}_1)$ is contained in a
connected component of the fixed point set for the action
of $T_1$ on $M$, and hence
${\mu}_T(\pi(\hat{M}_1))$ is contained in a wall for
the action of $T$ on $M$ with moment
map $\mu_T$. Recall that we have chosen $\xi \in \liets$ (respectively $\hat{\xi} \in \liets$)
to lie in a connected component of the set of regular values of $\mu_T$ (respectively
$\hat{\mu}_T$) containing 0 in its closure. It follows that for $k \muchgreater 0$ it is possible to reach
$\xi$ from $\hat{\xi}$ by crossing only those walls $\hat{\mu}_T(\hat{M}_1)$ for the
action of $T$ on $\hat{M}$ for which ${\mu}_T(\pi(\hat{M}_1))$ contains 0, and in particular
the constant value taken by $\mu_{T_1} \circ \pi$ is 0. Moreover we can always assume that the $\xi_1$
which appears in (\ref{new8}) lies in a connected component of the set of regular values of
$\hat{\mu}_{T/T_1}$ whose closure contains 0.

Suppose $\pi(\hat{M}_1)\cap M^{ss}= \emptyset$.
Since $T_1$ fixes $\hat{M}_1$, its image $\pi(\hat{M}_1)$ is contained in
a $T_1$-fixed point component in $M$. If the wall given by this component
does not pass through $0$, we do not need to cross the wall
determined by  $\hat{M}_1$
because it is a wall far away from $0$ in our scale.

Let $B\subseteq \mathbf{t}=\mathrm{Lie}(T)$ be a small ball
containing $0$. On $\mu^{-1}_T(B)$, we may choose an equivariant
differential form representing $\mathcal{D}|_{\mu^{-1}_T(B)}$,
supported in $\mu^{-1}_T(B)\cap M^{ss}$. If $\pi(\hat{M}_1)\cap
M^{ss}= \emptyset$ and the wall for $\pi(\hat{M}_1)$ passes
through $0$, we may assume $\xi_1 \in B$ and the wall crossing
term \beq \label{star} \kappa^{T/T_1}_{\hat{M}_1,\xi_1}( \res
_{\hat{M}_1}(\a\b \cald^2)) [ \hat{M}_1 /\!/_{\xi_1} (T/T_1)_c]
\eeq vanishes since $\alpha \beta
\mathcal{D}^2|_{\pi(\hat{M}_1)\cap \mu_T^{-1}(\xi_1) }=0$.

If $\pi(\hat{M}_1)\cap M^{ss}$ is non-empty but not contained in $GZ^{ss}_R$,
then the wall $\hat{\mu}_T(\hat{M_1})$
is determined by the  image of
$\hat{M}_1 \setminus \pi^{-1}(\overline{GZ_R^{ss}})$
under the moment map $\hat{\mu}_T$ which is determined by the linearization of the action on $\hat{M}$ given by the induced action on the
line bundle $\pi^* L^{\otimes k} \otimes {\mathcal O}(-E)$. As
above we can assume that $0 \in \mu_T(\pi(\hat{M_1})) $ and hence the constant
value taken on $\pi(\hat{M_1})$  by the $T_1$-moment map $\mu_{T_1}$
is $0$. This means that the induced action of $T_1$ on the restriction of
$L$ to $\pi(\hat{M_1})$ is trivial so the same is true for the
induced action of
$T_1$ on the restriction of $\pi^* L^{\otimes k} $ to
$\hat{M_1}$. Moreover since
$\hat{M_1} \setminus \pi^{-1} (\overline{GZ_R^{ss}})$ does not meet the exceptional
divisor $E$ we have
$${\mathcal O}(-E)|_{\hat{M_1} \setminus \pi^{-1} (\overline{GZ_R^{ss}} )}
\cong  {\mathcal O}|_{\hat{M_1} \setminus \pi^{-1} (\overline{GZ_R^{ss}} )}$$
and the induced action of $T_1$ on this is trivial since $T_1$
acts trivially on $\hat{M_1}$. Hence the induced action of
$T_1$ on the restriction of $\pi^* L^{\otimes k}  \otimes {\mathcal O}(-E)$
to $\hat{M_1} \setminus \pi^{-1}(\overline{GZ_R^{ss}}) $ is trivial,
and hence the constant value taken by $\hat{\mu}_{T_1} $ on
$\hat{M_1}$ is $0$. Therefore the wall $\hat{\mu}_{T}(\hat{M_1})$
passes through $0$, and hence need not be crossed.

Hence the only wall crossing terms we need to consider are for the
fixed point sets
$\hat{M}_1$ such that $\emptyset\neq
\pi(\hat{M}_1)\cap M^{ss}\subseteq GZ^{ss}_R$
as required. 

\begin{example}
Let us consider again the action of $G=SL(2)$ on $M=\PP_n$
when $n$ is even. We need to
blow $M$ up along the closure of the $G$-orbit of the unique
$T$-fixed point $p_0$ in $\mu^{-1}(0)$ in order to
obtain $\tilde{M}$. The only $T$-fixed point in $Gp_0$ is
$p_0$ itself and the weights of the induced action on the fiber
$\PP_{n-3}$ over $p_0$ in the exceptional
divisor of $\tilde{M}$ are $-n,2-n,...,n-2,n$ with $0, -2, 2$
omitted.
If as in Example \ref{caps} $\eta=q(\xi,\zeta^2) \in \hk (M)$ has degree $2(n-3)$ then
$$\k_{\tilde{M}}(\eta)[\xtg]=-\frac12
\k_{\tilde{M}}(\eta\mathcal{D}^2)[\tilde{M}
\git T_c]$$
and we can write
$$\k_{\tilde{M}}(\eta\mathcal{D}^2)[\tilde{M}
\git T_c]=\big(\k_{\tilde{M}}(\eta\mathcal{D}^2)[\tilde{M}
\git T_c]-\k_{\tilde{M}}(\eta\mathcal{D}^2)[\tilde{M}
\git_{\xi} T_c]\big) +\k_{\tilde{M}}(\eta\mathcal{D}^2)[\tilde{M}
\git_{\xi} T_c]$$
where $\xi$ is a positive number between $0$ and $1$.

Notice that $\k_{\tilde{M}}(\eta\mathcal{D}^2)[\tilde{M}
\git_{\xi} T_c]=\k_{M}(\eta\mathcal{D}^2)[{M}
\git_{\xi} T_c]$, which is given by summing up the wall crossing terms
for the components of the fixed point set with positive moment map values. Hence,
$$\k_{\tilde{M}}(\eta\mathcal{D}^2)[\tilde{M}
\git_{\xi} T_c]=-2\res_{X=0}\bigl( \sum_{j:0<n-2j\leq n}
\frac{q((n-2j)X,X^2)}{(2X)^{n-2} \prod_{k\neq j}(k-j)}\bigr).$$

On the other hand, the difference
$\k_{\tilde{M}}(\eta\mathcal{D}^2)[\tilde{M}
\git_0 T_c]-\k_{\tilde{M}}(\eta\mathcal{D}^2)[\tilde{M}
\git_{\xi} T_c]$ is given by the wall crossing terms for the components of the fixed point set
over $p_0$ with positive moment map values.
Hence the difference is
$$2\res_{X=0}\bigl(\sum_{j:2<n-2j\leq n}
\frac{q(0,X^2)}{2^{n-3} X^{n-2} (n-2j) \prod_{k\neq j,0}(k-j)}\bigr)$$
Note however that because $n$ is even the residue is always $0$.
Therefore, we get
$$\k_{\tilde{M}}(\eta)[\xtg]=\res_{X=0}\bigl( \sum_{j:0<n-2j\leq n}
\frac{q((n-2j)X,X^2)}{(2X)^{n-2} \prod_{k\neq j}(k-j)}\bigr).$$
\end{example}

\bigskip

Let us consider in more detail the walls $\hat{\mu}_T (\hat{M_1})$
and the wall-crossing terms (\ref{star}) in the general case.
 Consider the subset
$$L = \{g \in G: T_1 \subseteq  gRg^{-1} \}$$
of $G$.
By (8.10) of \cite{Ki2}, $L$ is the disjoint union of finitely many
double cosets for the left $N^{T_1}_0$ action and the right $N^R$
action
$$L=\bigsqcup_{1\le i\le m} N^{T_1}_0g_iN^R$$
where $N^{T_1}_0$ is the identity component of the normalizer
of $T_1$  in $G$ and $N^R$ is the normalizer of $R$ in $G$.
Since $G Z_R^{ss}  \cong G \times_{N^R} Z_R^{ss}$ it follows that
 the $T_1$-fixed point set in
$GZ^{ss}_R$ is the disjoint union
$$\bigsqcup_{1\le i\le m}N^{T_1}_0Z^{ss}_{R_i}$$
where $R_i=g_iRg_i^{-1}\supset T_1$.

\begin{definition}
Let $W_{i,j}\to N^{T_1}_0Z^{ss}_{R_i}$
for $1 \le j \le l_i$ be the $T_1$-eigenbundle
 of the restriction  to $N_0^{T_1} Z_{R_i}^{ss}$
of the normal bundle to  $GZ^{ss}_R$ on which $T_1$ acts with
 weight $\beta_{i,j}$.
\end{definition}

 Then $\bigsqcup_{i,j} \PP W_{i,j}$ is the  $T_1$ fixed point
set in $\pi^{-1}(GZ^{ss}_R)$, and we have proved

\begin{lemma} \label{l8.7}
Let $\hat{M_1}$ be a component of a fixed point set of a circle
subgroup $T_1$ satisfying
$ \emptyset \ne \pi(\hat{M_1}) \cap M^{ss} \subseteq G Z_R^{ss}$
as in Lemma \ref{l8.4}. Then
$$ \hat{M_1} \cap \pi^{-1} (M^{ss}) = \PP W_{i,j} $$
for some $i,j,$ and the corresponding wall crossing term
(\ref{star}) is
$$
\int_{\PP W_{i,j}\quott_{\xi_1}  {(T/T_1)_c} }
 \kappa^{T/T_1}_{\PP W_{i,j}, \xi_1}
\left ( {\rm res}_{X_1 = 0}  \frac{\alpha \beta \cald^2} {e_{\PP W_{i,j}}  }
 \right ) $$
where $e_{\PP W_{i,j}}$ is the $T$-equivariant Euler class of the normal bundle
to $\PP W_{i,j}$  in $\hat{M}$.
\end{lemma}

\begin{rem} \label{r8.8}
$\hat{M}_1$ is fixed by $T_1$, so its
images under $\hat{\mu}$ and $\mu \circ \pi$ are fixed by
$T_1$, and hence they are contained in the Lie algebra of the
centralizer of $T_1$. As $T_1$ is a circle it has no nontrivial
orientation preserving automorphisms, so the connected
component of this centralizer is just $N_0^{T_1}\cap K$.
In fact when $K$ is connected any such centralizer is connected, so the
centralizer of $T_1$ in $K$ is $N_0^{T_1}\cap K$ and the
centralizer of $T_1$ in $G$ is $N_0^{T_1}$.
Also the $T_1$ component of $\hat{\mu}$ is constant on $\hat{M}_1$;
let us call it $\xi_2$. Then it follows that
$$\hat{\mu}^{-1}_{N_0^{T_1}/T_1}(\xi_1) \cap \hat{M}_1 = \hat{\mu}^{-1}(\xi_1 + \xi_2) \cap \hat{M}_1.$$
Note that $\xi_1$ and $\xi_2$ can be taken to be arbitrarily close to 0, so that $\hat{\mu}^{-1}(\xi_1 + \xi_2)$
lies over semistable points in $M$ (cf. Remark 8.3).
\end{rem}

We are looking at the evaluation  of ${\rm res}_{\hat{M_1}} \alpha
\beta {\cal D}^2$ on
$$\left ( \hat{\mu}^{-1}_{T/T_1}(\xi_1) \cap \hat{M}_1\right  ) /(T/T_1).$$
This can be calculated inductively as follows. Let
$${\widetilde{ \PP  W}}_{i,j}    \quott N_0^{T_1} =
 {\widetilde{ \PP  W}}_{i,j}    \quott \left ( N_0^{T_1}/(T_1)_c\right )$$
be the partial desingularization of $\PP   W_{i,j}   \quott
N_0^{T_1}.$ By (\ref{new2}), we have
$$ \int_{{\widetilde{ \PP
W}}_{i,j}\quott_{\xi_1} (T/T_1)_c } \kappa^{T/T_1}_{{\widetilde{ \PP
W}}_{i,j} ,\xi_1 } \left ( {\rm res}_{X_1 = 0 } \frac{\alpha \beta
\cald^2}{e_{\PP  W_{i,j} } } \right ) =\int_{{ \PP
W_{i,j}}\quott_{\xi_1} (T/T_1)_c } \kappa^{T/T_1}_{{ \PP  W_{i,j}
},\xi_1 } \left ( {\rm res}_{X_1 = 0 } \frac{\alpha \beta
\cald^2}{e_{\PP  W_{i,j} } } \right ). $$

Since $\dim \hat{M}_1 < \dim M$ and $\dim \left (
N^{T_1}_0/(T_1)_c\right )<\dim G$, we can calculate the difference
between the wall crossing term above and \beq  \label{doublestar}
\int_{{\widetilde{ \PP W}}_{i,j}\quott (T/T_1)_c }
\kappa^{T/T_1}_{{\widetilde{ \PP W}}_{i,j}  } \left ( {\rm res}_{X_1
= 0 } \frac{\alpha \beta \cald^2}{e_{\PP W_{i,j} } } \right ) \eeq
inductively using the same procedure.

Now, as at (8.1), we find that (\ref{doublestar}) is equal up to a
scalar factor analogous to $\frac{n_0 (-1)^{n_+}}{n_0^{T} |W|}$ to
\begin{equation} \label{dagger}
\int_{{\widetilde{ \PP  W}}_{i,j} \quott N_0^{T_1} }
\kappa_{{\widetilde{ \PP W}}_{i,j}  }  \left ( {\rm res}_{X_1 = 0
} \frac{\alpha \beta \cald^2} { e_{\PP W_{i,j}}
(\cald_{N_0^{T_1}})^2 } \right )
\end{equation}
where
$\cald_{N_0^{T_1}}$ is the product of the positive
roots of $N_0^{T_1}$.

\begin{rem}\label{extra2}
The wall crossing term (\ref{star}) is an integral over the
quotient $\hat{M}_1/\!/_{\xi_1} (T/T_1)_c = \hat{\mu}_T^{-1}(\xi_1 + \xi_2) \cap \hat{M}_1 /T$.
Since the intersection of $\hat{M}_1$ with $\pi^{-1}(M^{ss})$ is a projective
bundle over $N_0^{T_1}Z_{R_1}^{ss}$, it is easier to describe
quotients of $\hat{M_1}$ by $N_0^{T_1}$ than by $T$, and we can use (1.6) to
relate integrals over quotients of $\hat{M}_1$ by $T$ to integrals over quotients of
$\hat{M}_1$ by $N_0^{T_1}$. Unfortunately, however, the quotient
$$\hat{M}_1 /\!/_{\xi_1} (N_0^{T_1}/(T_1)_c) = \hat{\mu}^{-1}(\xi_1 + \xi_2) \cap \hat{M}_1/(N_0^{T_1} \cap K)$$
is only well defined if $\xi_1$ is centralized by $N_0^{T_1}$, and this is not necessarily the case.
To overcome this problem, we can use induction to relate our integral over $\hat{M}_1/\!/_{\xi_1} (T/T_1)_c $
to an integral over the partial desingularization of $\hat{M}_1/\!/ N_0^{T_1}$, as above, or indeed to an
integral over the partial desingularization of $\hat{M}_1/\!/_{\zeta} N_0^{T_1}$ for any $\zeta$  in the Lie
algebra of the centre of $N_0^{T_1} \cap K$. In what follows we shall assume for simplicity that $\zeta =0$,
but the same argument will work for any central $\zeta$. In particular if $\xi_1$ itself is central then we can
simply work with $\hat{M}_1/\!/_{\xi_1} N_0^{T_1}$ and no partial desingularization is required.
\end{rem}

\begin{rem} \label{r8.9}
At this point it is useful to recall from \cite{Ki2}, Lemma 7.8 that
the fibre at $x \in Z_R^{ss}$ of the restriction
$\pi: \hat{M}^{ss} \to M^{ss}$ of $\pi $
to $\hat{M}^{ss}$ is the semistable  stratum of the projective space
associated to the normal to $G Z_R^{ss}$ in $M^{ss}$ at $x$
with respect to the representation of $R$ on the normal
space. Since $R \subseteq N^R \subset G$  and semistability for some
reductive group implies semistability for any reductive subgroup, it follows that in such a fibre
semistability for $R$, for $N^R$ and for $G$ are all equivalent. Similarly
(cf. Remark \ref{r8.8}) in the intersection
$ \PP ( W_{i,j})_x$ of $\hat{M_1} =
 \PP  W_{i,j}$ with the fibre
of $\pi$ at $x \in Z^{ss}_{R_i} $, semistability for $R_i \cap
N_0^{T_1}$, for $N^{R_i}  \cap N_0^{T_1}$ and for $ N_0^{T_1}$ are
all equivalent.
\end{rem}

Recall that $G Z_R^{ss} \cong G \times_{N^R} Z_R^{ss}$
and hence
$$N_0^{T_1} Z_{R_i}^{ss} \cong N_0^{T_1} \times_{N_0^{T_1} \cap N^{R_i}}
Z_{R_i}^{ss}$$
and so
$$  \PP  W_{i,j} \cong
 N_0^{T_1} \times_{N_0^{T_1} \cap N^{R_i}}
\PP  W_{i,j}\mid_{Z_{R_i}^{ss} }     $$
$$   {\widetilde{\PP  W}}_{i,j} \cong
 N_0^{T_1} \times_{N_0^{T_1} \cap N^{R_i}}
{\widetilde{\PP  W}}_{i,j}\mid_{Z_{R_i}^{ss} }. $$ It then follows
immediately from Remark \ref{r8.9} that \beq \label{triplestar}
 {\widetilde{\PP  W}}_{i,j}\quott N_0^{T_1} \cong
{\widetilde{\PP  W}}_{i,j}\mid_{Z_{R_i}^{ss} } \quott N_0^{T_1} \cap
N^{R_i}. \eeq
\begin{lemma} \label{l8.10}
$R_i (N_0^{T_1} \cap N^{R_i})$ is a subgroup of finite index in $N^{R_i}$.
\end{lemma}
\Proof We modify the proof of \cite{Ki2}, Lemma 8.10. First note that each
subgroup $R_i$ and $N_0^{T_1}$ and $N^{R_i}$ of $G$ is the complexification of
its intersection with the maximal compact
subgroup $K$ of $G$, and so it suffices to show that
$(R_i \cap K) (N_0^{T_1} \cap N^{R_i} \cap K) $ has finite
index in $N^{R_i} \cap K$.

Consider the action of $R_i \cap K$ on the homogeneous space
$(N^{R_i} \cap K)/T_1$. If $k \in N^{R_i } \cap K$ then the stabilizer
in $R_i \cap K$ of the coset $k T_1$ is
$$ \{ s \in R_i \cap K \mid sk T_1 = k T_1 \} = R_i \cap k T_1 k^{-1}$$
since $T_1 \subseteq R_i$ and $k$ normalizes $T_1$. But only finitely
many conjugacy classes of subgroups of $R_i \cap K$ can occur as stabilizers
of elements of $(N^{R_i} \cap K)/T_1 $  (see \cite{Mo}). Therefore there
exist $k_1, \dots, k_m \in N^{R_i} \cap K$ such that if
$k \in N^{R_i } \cap K $ then
$$ k T_1 k^{-1} = r k_j T_1 k_j^{-1} r^{-1}$$
for some $r \in R_i \cap K$ and $j \in \{ 1, \dots, n\}$,
and hence
$$k \in (R_i \cap K) k_j (N^{T_1}\cap N^{R_i} \cap K) $$
$$ = k_j (R_i \cap K) (N^{T_1} \cap N^{R_i} \cap K) $$
since $k_j$ normalizes $R_i \cap K$. Since $N^{T_1}_0$ has finite
index in $N^{T_1}$, the result follows.

\begin{rem} \label{r8.11} It follows immediately that
the connected component $N_0^{R_i}$ of $N^{R_i}$ is equal to the connected
component $R_i (N_0^{T_1} \cap N^{R_i} )_0$ of $R_i (N_0^{T_1} \cap N^{R_i} )$.
\end{rem}

\begin{corollary} \label{c8.12}
There is a natural isomorphism
$$ Z_{R_i} \quott (N_0^{T_1} \cap N^{R_i} )_0 \cong Z_{R_i} \quott N_0^{R_i} $$
and  finite-to-one surjections
$$ Z_{R_i} \quott N_0^{R_i} \to  Z_{R_i} \quott N^{R_i} $$
and
$$  Z_{R_i} \quott N_0^{R_i} \to Z_{R_i} \quott N_0^{T_1} \cap N^{R_i}. $$
\end{corollary}
\Proof This follows immediately because the linear action of
$R_i$ on $Z_{R_i}$
is trivial.

\bigskip

We have reduced the calculation of the wall crossing term
(\ref{star}) to the calculation of the integral (\ref{dagger})
over $  {\widetilde{\PP W}}_{i,j} \quott N_0^{T_1} $ and from
(\ref{triplestar}) we know that
$$  {\widetilde{\PP W}}_{i,j} \quott N_0^{T_1}
\cong  {\widetilde{\PP W}}_{i,j} \mid_{ Z_{R_i}^{ss}}
 \quott N_0^{T_1} \cap N^{R_i}.$$

We have a surjection
\begin{equation} \label{commutdiag}
 {\widetilde{\PP W}}_{i,j}  \mid_{Z_{R_i}^{ss} } \quott N_0^{T_1} \cap N^{R_i}
\stackrel{\Psi}{\longrightarrow} Z_{R_i} \quott N_0^{T_1} \cap N^{R_i}.
\end{equation}
\begin{lemma} \label{l8.13}
The fibre of $\Psi$ is the partial desingularization
$${\widetilde{ \PP W}}_{i,j}|_x \quott Stab (x)\cap N_0^{T_1} \cap
N^{R_i} $$ of the quotient of the projective space $\PP
W_{i,j}|_x$ by $Stab (x)\cap N_0^{T_1} \cap N^{R_i} $.
\end{lemma}
\Proof This follows since the centre of the blowup required at any
stage of the construction of ${\widetilde{\PP W}}_{i,j} \quott
N_0^{T_1}$ intersects the fibre at $x \in Z^{ss}_{R_i} $ in the
centre of the corresponding blowup required for the construction
of the partial desingularization
$$ {\widetilde{ \PP W}_{i,j} } |_x  \quott Stab
(x)\cap N_0^{T_1} \cap N^{R_i} $$ (cf. the proofs of \cite{Ki3}
Lemma 1.16 and Proposition 1.20).

We can now perform the integration needed over
$$ {\widetilde{\PP W}}_{i,j}|_{Z^{ss}_{R_i}} \quott N_0^{T_1} \cap N^{R_i} $$
by first integrating over the fibres of $\Psi $. This involves
integrals over partial desingularizations
$$ {\widetilde{ \PP W}_{i,j} }|_x
 \quott Stab(x) \cap N_0^{T_1} \cap N^{R_i} $$
and
$$ Z_{R_i} \quott N_0^{T_1} \cap N^{R_i}$$
which can be calculated inductively.
\begin{rem} \label{r8.14}
(i) Recall  from Corollary \ref{c8.12}
that there is a natural surjection $ Z_{R_i} \quott N_0^{R_i}
\stackrel{\Phi}{\longrightarrow}   Z_{R_i} \quott N^{R_i}  \cap N_0^{T_1} $
such that if $ \eta \in H^*(Z_{R_i} \quott N_0^{T_1} \cap N^{R_i})$
then
$$ \int_{Z_{R_i} \quott N^{R_i}  \cap N_0^{T_1} } \eta
= \frac{1}{c} \int_{Z_{R_i} \quott N_0^{R_i} }
\Phi^* \eta, $$
where $c$ is the order of
$$\frac{Stab_{R_i(N_0^{T_1 } \cap N^{R_i}) }(x)}{Stab_{N_0^{R_i}} (x)}$$
for  a generic $x \in Z^{ss}_{R_i}$.

(ii) The connected component of $Stab (x)$ for $x \in
Z_{R_i}^{ss}$ is $R_i$ so $N_0^{T_1} \cap R_i$ is a subgroup of
finite index in $Stab(x) \cap N_0^{T_1} \cap N^{R_i},$ and hence
as in (i) we can reduce integrals over
$$
{\widetilde{ \PP W}_{i,j} }|_x \quott Stab (x) \cap N_0^{T_1} \cap
N^{R_i} $$ to integrals over
$$
{\widetilde{\PP W}_{i,j}}|_x \quott R_i \cap N_0^{T_1}. $$

(iii) We have assumed for simplicity that the blowup along $G
Z_{R_i}^{ss}$ is the first stage of the construction of $\tilde{M}
\quott G$ from $M \quott G$. If instead we are at a later stage of
the procedure, we need to replace $Z_{R_i} \quott N^{R_i} $ by its
partial desingularization $\tilde{Z}_{R_i} \quott N^{R_i} $
throughout; cf. \cite{Ki2}, Lemma 8.8.
\end{rem}

Recall that in order to calculate the pairing $\kappa_{\tilde{M}}(\alpha \beta)[\tilde{M}/\!/G]$, by
(8.1) and Lemma 8.2 it suffices to calculate the difference
\begin{equation} \label{dbdagger} \kappa^T_{\tilde{M},\xi}
(\alpha \beta {\cal D}^2)[\tilde{M}/\!/_\xi T_c] -
\kappa^T_{\tilde{M}} (\alpha \beta {\cal D}^2)[\tilde{M}/\!/ T_c]
\end{equation}
which is a sum of differences of the form
$$ \kappa^T_{\hat{M},\xi} (\alpha \beta {\cal D}^2)[\hat{M}/\!/_\xi T_c]
- \kappa^T_{\hat{M},\hat{\xi}} (\alpha \beta {\cal D}^2)[\hat{M}/\!/_{\hat{\xi}} T_c] $$
fo suitable $\xi$ and $\hat{\xi}$. Recall also that these differences are in turn sums
of wall crossing terms (\ref{star})
$$\kappa^{T/T_1}_{\hat{M}_1,\xi_1}( \res _{\hat{M}_1}(\a\b \cald^2)) [ \hat{M}_1 /\!/_{\xi_1} (T/T_1)_c]
=
\int_{\PP W_{i,j}\quott_{\xi_1}  {(T/T_1)_c} }
 \kappa^{T/T_1}_{\PP W_{i,j}, \xi_1}
\left ( {\rm res}_{X_1 = 0}  \frac{\alpha \beta \cald^2} {e_{\PP W_{i,j}}  }
 \right ), $$
and the difference between any such wall crossing term and the corresponding integral (\ref{dagger})
$$
\int_{{\widetilde{ \PP  W}}_{i,j} \quott N_0^{T_1} }
\kappa_{{\widetilde{ \PP W}}_{i,j}  }  \left ( {\rm res}_{X_1 = 0
} \frac{\alpha \beta \cald^2} { e_{\PP W_{i,j}}
(\cald_{N_0^{T_1}})^2 } \right )
$$
over the partial desingularization of ${ \PP  W}_{i,j} \quott N_0^{T_1}$ can be calculated inductively.
Finally we have just observed that we can perform this last integral by integrating over the fibres of
the map
$$\Psi:
 {\widetilde{\PP W}}_{i,j}   \quott N_0^{T_1}
\rightarrow Z_{R_i} \quott N_0^{T_1} \cap N^{R_i},$$
and the integral over the fibre of $\Psi$ at a point represented by $x \in Z_{R_i}^{ss}$ reduces to an
integral over the partial desingularization ${\widetilde{\PP W}}_{i,j}
  |_x \quott R_i$ of the GIT quotient
of the projective space ${\PP W}_{i,j}   |_x$ by $R_i$, while the
resulting integral over the base space $ Z_{R_i} \quott N_0^{T_1}
\cap N^{R_i}$ can be pulled back to an integral over $ Z_{R_i}
\quott  N^{R_i}$. We can use induction and the residue formula
Theorem 3.1(i) to calculate integrals over $ Z_{R_i} \quott
N^{R_i}$ of any cohomology classes represented by equivariant
cohomology classes on $Z_{R_i}$, so we can calculate
(\ref{dagger}) provided that we are able to represent the
cohomology class on $ Z_{R_i} \quott  N^{R_i}$ obtained by
integrating
\begin{equation} \label{kappa}
\kappa_{{\widetilde{ \PP W}}_{i,j}  }  \left ( {\rm res}_{X_1 = 0
} \frac{\alpha \beta \cald^2} { e_{\PP W_{i,j}}
(\cald_{N_0^{T_1}})^2 }\right )
\end{equation}
over the fibres of $\Psi$ by an equivariant cohomology class on $Z_{R_i}$. But if $x \in Z_{R_i}^{ss}$
then we can use induction and the residue formula Theorem 3.1(i) to express the integral over
${\widetilde{ \PP W}}_{i,j}  |_x \quott R_i$ of
(\ref{kappa}) 
in terms of residues of integrals over fibres at $x$ of projective subbundles of
$\PP W_{i,j} \mid _{Z_{R_i}^{ss}}$ which extend naturally to projective bundles over
$Z_{R_i}$. Such integrals can be represented by equivariant cohomology classes
on $Z_{R_i}$ using the following lemma, which therefore completes our procedure
for calculating $\kappa_{\tilde{M}}(\alpha \beta)[\tilde{M}/\!/G]$.

\begin{lemma}
 Let $E$ be a rank $r$ complex vector bundle over
a manifold $M$. Then
$$ \int_{{\bf P} E} \eta = {\rm Res}_{y = 0 } \int_M \frac{\eta}{p(y)}
$$
where $p(y) =  y^r +
c_1(E) y^{r-1} + \dots + c_r(E).  $
\end{lemma}

{\bf Proof:} This follows from standard arguments in algebraic topology.
Here we give a simple argument. Let $P$ denote ${\bf P}(E)$, and let
$\pi: P \to M$ be the projection.
Then $H^*(P)$ is isomorphic to the quotient of the polynomial ring $ H^*(M) [y]$
by the ideal generated by $p(y)$ as above. For any
cohomology class $\eta \in H^*(P)$, we have
$\eta = \sum_{i=0}^{r-1} \pi^* \beta_i y^i$ for suitable classes $\beta_i \in
H^*(M)$, where we have identified $y$ with the first Chern class of the
hyperplane line bundle $L \to P$.
Then
$$ 
 \frac{1}{p(y)} = \frac{1}{y^r(1 + \frac{c_1(E)}{y} + \dots +
\frac{c_r(E)}{y^r} )} = \sum_{k \geq 0} \frac{(-1)^k}{y^r}\left(\frac{c_1(E)}{y} + \dots + \frac{c_r(E)}{y^r} \right)^k $$
and so
$$\int_P \eta = \int_M \pi^* \beta_{n-1} = 
{\rm Res}_{y = 0 } \int_M \frac{\eta}{p(y)}.
$$

\renorm
\section{Witten's integral}
In the case when $0$ is a regular value of the moment map
$\mu$ Witten relates the intersection pairings of two classes 
$\k_M(\a)$, $\k_M(\b)$ of
complementary degrees in $H^*(M/\!/G)$
coming from $\a$, $\b $ $\in H^*_K(M)$ to the asymptotic behaviour 
of the integral $\cali^\epsilon (\a \b e^{i\baromega})   $ given by 
\beq \label{iepstwo} \cali^\epsilon (\eta e^{i\baromega}) 
=\frac{1}{(2 \pi )^s  \vol(K) }  \int_{X \in \liek} [d X] e^{-\e<X,X>/2}
\int_M \eta(X) e^{i \omega} e^{i \mu(X)},  \end{equation}
where as before $\baromega = \omega + \mu$. He expresses the 
integral as a sum of local contributions, one of which reduces to
the intersection pairing required while the rest tend
to 0 exponentially fast as $\e$ tends to 0.

Even when $0$ is not a regular value of $\mu$, Witten's integral 
$\cali^\epsilon (\eta e^{i\baromega}) $
decomposes into the sum of a term $\cali^\epsilon_0(\eta 
e^{i\baromega})$
determined by the action of $K$ on an arbitrarily small 
neighbourhood
of $\mu^{-1}(0)$, and other terms which tend to zero exponentially
fast as $\epsilon \to 0 $. We shall see that there is a residue formula 
for $\cali^\epsilon_0(\eta e^{i\baromega})$ which again is a sum over 
components of the fixed point set of $T$ on $M$. This residue 
formula is related to the formulas for pairings in the intersection homology
of $M/\!/G$ given in previous sections,  but it
is not in general a polynomial
in $\epsilon$; instead it is a polynomial in $\sqrt{\e}$, as has been proved by
Paradan (see \cite{paradan} Cor.5.2).

In Sections 4  and 7 of \cite{JK1} it is proved
that the integral $\cali^\epsilon (\eta e^{i\baromega})$ can be 
expressed
as 
\beq \label{idef} \cali^\epsilon (\eta e^{i\baromega}) =  
\frac{1}{(2 \pi )^l \epsilon^{s/2}|W| \vol(T) }\int_{y \in \bf{t}^*} dy 
e^{-<y,y>/{2 \epsilon} } Q^\eta(y)
\end{equation}
where $W$ is the Weyl group of $K$ and $Q^\eta(\cdot)$ is a 
piecewise polynomial function on ${\bf t}^*$ supported on cones 
each of which has its apex at $\mu_T(F)$ for some component 
$F$ of the fixed point 
set of $T$ on $M$. Here as before $s$ is the dimension of $K$, while $l$ 
is the dimension of the maximal torus $T$.
For the definition of $Q^\eta$ see 
the statement of Theorem 7.1 of \cite{JK1}: it is 
$$Q^\eta = 
\cald F_T(\cald  \sigma^\eta) $$
where 
$$ \sigma^\eta(X)  = 
\Pi_* (e^{i \omega } e^{i \mu(X)} \eta(X)),  $$
and $\Pi_*$ denotes the integral over $M$ while $F_T $ is the 
Fourier transform over ${\bf t}$. Equivalently 
if $\{e_j\}$ is a basis for ${\bf t}$
so that an element of ${\bf t}$ is given as
$y = \sum_j y_j e_j $,
 we may write 
$$Q^\eta(y) = \cald \calp_{\cald} (F_T  \Pi_* \sigma^\eta)$$ 
where we define the differential operator
$$\calp_\cald = \prod_{\gamma> 0 } \sum_j \gamma (e_j) \frac{\partial}{\partial y_j} $$
as a product over the positive roots $\gamma$ of $K$.

We shall  need to introduce a set  of (possibly degenerate) cones $\calc = \{ C_1,
 \dots, C_d \}, $ each with apex at $0$, for which $\liets$ is the 
union of $C \in {\mathcal{C}}$,   the intersection 
of any two is contained in their boundaries and 
$Q^\eta$ is polynomial on a neighbourhood of $0 $ in each $C$.
Let $Q^\eta_0$ be the piecewise polynomial function which is 
polynomial on each cone 
$C \in \calc$ (each
having its apex at $0$) and which coincides with $Q^\eta$ near 
$0$. 
Let \beq \label{9.1} {\cali^\epsilon}_0 (\eta 
e^{i\baromega}) =   \frac{1}{ (2 \pi )^l \epsilon^{s/2 } |W| \vol(T) } 
 \int_{y \in {\bf t}^* } [dy] e^{-<y,y>/{2 \epsilon} } Q^\eta_0(y).
\end{equation}
Then the argument of \cite{JK1} Section 6 shows that there exist real 
numbers $\rho_\beta>0$ and functions $h_\beta: \RR^+ \to \RR$ such
 that for some $N_\beta \ge 0 $ the product $\epsilon^{N_\beta} 
h_\beta (\epsilon)$ remains bounded as $\epsilon \to 0^+$ and 
\beq |\cali^\epsilon (\eta e^{i\baromega}) - \cali_0^\epsilon(\eta 
e^{i\baromega})| \le \sum_\beta e^{- \rho_\beta/{2 \epsilon} }
h_\beta(\epsilon).  \end{equation}

In Section 8 of \cite{JK1} a residue formula is given for 
$\cali^\epsilon_0 (\eta e^{i\baromega})$ in the case when $0$
is a regular value of the moment map $\mu$. In this case 
$n_0 \cali^\epsilon_0 (\eta e^{i\baromega})$ is equal to 
$ \eta_0 e^{i\omega + \epsilon \Theta}[M/\!/K] $, 
where $\Theta\in H^*(\xg)$ was
defined at (1.3) and $n_0$ is the order of the stabilizer of a generic
point of $\mu^{-1}(0)$, and thus $I^\epsilon_0(\eta e^{i \bom})$ is a polynomial function of 
$\epsilon$.  The proof of Theorem 8.1 of \cite{JK1} can be modified 
to obtain a formula for $\cali^\epsilon_0 (\eta e^{i \baromega})$ in 
the general case  when it may not be a polynomial in $\epsilon$
(see Example 9.7 below). 

If $F \in \calf$ is a component of the fixed 
point set of the maximal torus
of $T$ acting on $M$ let $\beta_{F,j}$ (for $j\in J_F$) be the weights 
of the action of $T$ on the normal bundle to $F$ in $M$. We choose 
a connected component $\Lambda$ of the set of $\xi \in \liet$ for which
$\beta_{F,j}(\xi) \ne 0 $ for all $F$ and $j$; we then adjust the signs 
of the $\beta_{F,j}$ (for all $F$ and $j$) in such a way that 
$\beta_{F,j}(\xi) > 0$ for all $\xi \in \Lambda$ (see \cite{GLS}).

We can then define $\res^{\Lambda}(h(X)[dX])$ as in \cite{JK1} Section 8
when $h(X)$ is of the form
$$h_{\lambda}(X) = \frac{q(X) e^{i \l (X)}}{\beta_1(X)...\beta_N(X)},$$
where $q(X)$ is a polynomial in $X\in \liet$ while $\beta_1,...,\beta_N \in \liets$
all lie in the dual cone of $\Lambda$ and $\l \in \liets$ does not lie in 
any cone of dimension at most $l-1$ spanned by a subset of $\{\beta_1,
..., \beta_N\}$. By \cite{locquant} Proposition 3.2 it is uniquely determined
by the following properties:

\noindent i) If $\{ \beta_1,...,\beta_N \}$ does not span $\liets$ as a
vector space then $\res^{\Lambda}(h_{\lambda}(X)[dX])=0$.

\noindent ii) $\res^{\Lambda}(h_{\lambda}(X)[dX])$ $= \sum_{m\geq 0}
\lim_{s \to 0+} \res^{\Lambda}(\frac{(i\l (X))^m}{m!} h_{s\l}(X)[dX])$.

\noindent iii) If $q(X) = X_1^{j_1}...X_l^{j_l}$ then the limit
$\lim_{s \to 0+} \res^{\Lambda}(h_{s\l}(X)[dX])$ is 0 unless 
$N=l+j_1 + .... + j_l$.

\noindent iv) If $q(X) = 1$ and $N=l$ and $\{ \beta_1,...,\beta_l \}$ is
a basis for $\liets$ then $\res^{\Lambda}(h_{\l}(X)[dX])=0$ unless
$\l = \l_1 \beta_1 +...+ \l_l \beta_l$ where $\l_j >0$ for each $j$, and
if this is the case then $\res^{\Lambda}(h_{\l}(X)[dX])=|\det \bar{\beta}|^{-1}$,
where $\bar{\beta}$ is an $l \times l$ matrix whose columns are the coordinates
of $\beta_1,...,\beta_l$ with respect to any orthonormal basis of $\liet$.

Finally in order to remove the restriction on $\lambda$ we choose $\rho \in \liets$
such that $-\rho$ lies in the dual cone of $\Lambda$ and define
$$\res^{\rho,\Lambda}(h_{\l}(X)[dX]) = \lim_{s \to 0+} \res^{\Lambda}(h_{\l + s\rho}(X)[dX]).$$
Except for the additional factors of $i$ discussed in Remark 3.4, when applied to suitable 
meromorphic differential forms on the complexified Lie algebra of the torus, $\res^{\rho,\Lambda}$
gives the multivariable residue which appeared in Section 3.

If $F \in \calf$ we  define a cone $C(F)$ in $\bf{t}$,  with apex 
at $\mu_T(F)$, by 
\beq C(F) = \{\mu_T(F) - \sum_j s_j \beta_{F,j}: s_j \geq 0 \}  
\end{equation}
where the $\beta_{F,j}$ are the weights of the action of $T$ on the 
normal bundle to $F$ with adjusted signs as above.  By subdividing 
the cones $C \in \calc$ if necessary, we can assume that for each 
$F \in \calf$ and each $C \in \calc$, 
either there is a neighbourhood of 0 in $C$ which is contained in 
$C(F)$
or there is a neighbourhood of 0 in $C$ which does not meet the 
interior of $C(F)$. Then if $F \in \calf$ we let $\calc_F$ be the 
set of $C \in \calc$ such that there is a neighbourhood of 0 in $C$ 
which is contained in $C(F)$.

\newcommand{\inpr}[1]{ \langle #1 \rangle }

\begin{theorem} \label{t9.1}
If $\eta \in H^*_K(M)$ then $ \cali^\epsilon_0(\eta e^{i\baromega}) $ is equal to
$$ \frac{A_K}{\epsilon^{s/2}}  \sum_{F \in {\mathcal F} } \sum_{C \in {\mathcal C_F}}  
\int_{y \in C  } [dy]  \cald(y) e^{- <y, y>/{2\epsilon} } \calp_C 
\res^{\rho,\Lambda} \Bigl (     \cald(X) \int_F\frac{i_F^* (\eta(X)  
e^{i \omega}) }{e_F(X) } e^{i <\mu(F)-y,X>} \Bigr )  $$
where  the constant $A_K$ is given by
$$ A_K = \frac{i^l(2 \pi)^{-l/2}  }{|W|\vol(T)}. $$
Here $e_F$ denotes, as before, the equivariant
Euler class of the normal bundle to $F$, and 
if $f$ is a piecewise 
polynomial function which is polynomial on a neighbourhood
of $0$ in the cone $C$, then $\calp_C(f) $ denotes the polynomial which is 
equal to $f$ on a neighbourhood of $0$ in $C$. 
\end{theorem}

\newcommand{\hf}{H_{C,F}^{\rho,\Lambda,\eta}}
\newcommand{\hfo}{H_{C,F,0}^{\rho,\Lambda,\eta}}

\noindent{\bf{Proof:}} The proof is a 
straightforward modification of \cite{JK1}, Sections 4 and  8.
By (\ref{iepstwo}), we have 
$$I^\epsilon(\eta e^{i\bom}) = \frac{1}{(2 \pi)^s \vol K} \int_{X \in \liek}[dX ]
e^{- \epsilon |X|^2/2} \int_M \eta(X) e^{i \bom(X)}, $$
$$ = \frac{1}{(2 \pi)^l |W| \vol T} \int_{X \in \liet} [dX ]
e^{- \epsilon |X|^2/2} \cald(X)^2 \int_M \eta(X) e^{i \bom(X)}. $$
By Parseval's theorem this becomes 
$$  I^\epsilon(\eta e^{i\bom})  = \frac{1}{(2 \pi)^l \epsilon^{s/2} |W| \vol T} 
\int_{y \in \liet^*} [dy] \cald(y) e^{- |y|^2/{2 \epsilon}} F_T 
\Bigl (\cald(X) \int_M \eta(X) e^{i\omega}e^{i \mu(X) }\Bigr ), $$
where $F_T$ denotes the Fourier transform over $\liet$. 
We now expand the integral over $M$ using the abelian 
localization theorem, to get a
sum of terms
$$\sum_{F \in \calf} \int_F\frac{i_F^* (\eta(X)  e^{i \omega}) }{e_F(X) } e^{i \inpr{\mu(F),X}},$$
each corresponding to a component $F$ of the fixed point set
 of the action of $T$ (see for example \cite{GLS}). The Fourier transform can thus
be expressed as a sum over $\calf$ such that the term 
corresponding to $F \in \calf$ is a piecewise
polynomial function supported on a cone with apex at $\mu_T(F)$. 
Decomposing $\liets$ into cones $C$ with apex at $0$ as above, we find that 
the integrand on each cone $C$ is $e^{- |y|^2/{2 \epsilon} }$ times 
a piecewise polynomial function $p_C$. 
The functions 
$$\cald(y)\res^{\rho,\Lambda} \Bigl (   \cald(X) 
\int_F\frac{i_F^* (\eta(X)  e^{i \omega}) }{e_F(X) } e^{i \inpr{\mu(F)-y,X}}
 \Bigr)$$
are also piecewise polynomial, 
and we can assume that they 
are also polynomial on  a neighbourhood of $0$ 
in each of  the cones $C$ $\in \calc$. 
To define $I^\epsilon_0(\eta e^{i \bom})$, we replace  $p_C$ by the polynomial 
$p_C^0$ which equals $p_C$ in a neighbourhood of $0$ in $C$,
which by the argument of \cite{JK1} Section 8 is
\begin{equation}  \label{iintexpl}
(2 \pi)^{l/2} i^l				
\cald(y) 
\calp_C \res^{\rho,\Lambda} \Bigl (  
 \cald(X) \int_F\frac{i_F^* (\eta(X)  e^{i \omega}) }{e_F(X) } 
e^{i \inpr{\mu(F)-y,X} }  \Bigr). \end{equation}
This gives us the formula in the statement of the theorem.

\begin{rem} If the condition that $C \in \calc_F $ (in other words that a 
neighbourhood of $0$ in $C$ lies in $C(F)$) guarantees that the residue
$$\res^{\rho,\Lambda} \Bigl (   \cald(X) \int_F\frac{i_F^* (\eta(X)  
e^{i \omega}) }{e_F(X) } e^{i \inpr{\mu(F) -y,X} }  \Bigr) $$ 
is a polynomial function of $y$ on $C$ (not merely piecewise 
polynomial), then we can omit the expression $\calp_C$ from 
(\ref{iintexpl}).
\end{rem}

\begin{rem}
If $F \in \calf$ is such that $\mu_T(F)$
does not lie on a wall through 0 (or a wall
such that the affine hyperplane spanned by the wall passes
through 0), then 0 does not lie on the boundary of the cone $C(F)$, 
and hence either $\calc_F = \emptyset$ or $\calc_F = \calc$. 
If $\calc_F = \emptyset$ then $F \notin \calf_+$ and $F$ contributes 
zero to the expression for $ \cali^\epsilon_0(\eta e^{i\baromega})$ in 
Theorem 9.1. If $\calc_F = \calc$, so that $F \in \calf_+$, then the 
contribution of $F$ can be written as an integral over 
$\cup_{C \in \calc}C = \liet$ and by the arguments of \cite{JK1} Section
8 it is given by the same formula 
$$ C_K~ \res^{\rho,\Lambda}(\cald(X)^2\int_F \frac{i^*_F(\eta e^{i\bom})(X)
e^{- \epsilon<X,X>/2 } } {e_F(X)}[dX])$$
as in the case when $0$ is a regular value of $\mu$.
\end{rem}
\begin{example} \label{ex9.2} 
If $l = 1$ and there are two cones (both half lines) then when 
calculating $\cali^\epsilon_0 (\eta e^{i\baromega}) $ we need 
integrals of the form 
$$\int_{y \in C} e^{ - <y, y>/{2 \epsilon} }y^j dy  = \pm 
\int_0^\infty e^{ - y^2/{2 \epsilon}}  y^j dy. 
$$
Ignoring the sign, this can be evaluated as 
$$ \left [ - \epsilon e^{- y^2/{2 \epsilon} } y^{j-1}  \right]_0^\infty + 
\epsilon (j-1) \int_0^\infty e^{-y^2/{2 \epsilon} } y^{j-2} dy$$
$$ = \epsilon(j-1) \int_0^\infty e^{-y^2/{2 \epsilon} } y^{j-2} dy  $$
if $j>1$. By induction this equals 
$$ \epsilon^{(j- 1)/2} (j-1)(j-3) \dots 4 \cdot 2 \int_0^\infty 
e^{ - y^2/{2 \epsilon}}  y dy$$
$$ =\epsilon^{(j+ 1)/2} (j-1)(j-3) \dots 4 \cdot 2  $$
if $j$ is odd, and if $j$ is even it is 
$$\epsilon^{j/2} (j-1) (j-3) \dots 3 \cdot 1 \int_0^\infty 
e^{-y^2/{2\epsilon}} dy$$
$$ = \sqrt{\frac{\pi}{2}} \epsilon^{(j+1)/2} (j-1)(j-3) \dots 
3 \cdot 1   .        $$
Thus we expect that $\sqrt{\epsilon}$ will appear in the 
answer and $\cali^\epsilon_0(\eta e^{i\baromega}) $ will 
not in general be a polynomial in $\epsilon$. 
\end{example}

\begin{rem} By subdividing the cones $C \in \calc$ if 
necessary, we can assume that each $C \in \calc$ is of the form
$$C = \{s_1 b_1 +... + s_l b_l:s_1,...,s_l \in \RR, s_1 \ge 0,...,
s_m \ge 0\}$$
for some basis $b_1,...,b_l$ of $\liet$ and some $m \in \{0,...,l\}$.
Then the formula of Theorem 9.1 can be expressed as a linear 
combination 
(whose coefficients are independent of $\e$) of integrals of the form
$$\int_{y \in C} [dy] p(y) e^{-\langle y,y \rangle / 2 \e}$$
where $p(y)$ is a polynomial function of $y \in \liet$. Changing 
coordinates 
using the basis $b_1,...,b_l$ of $\liet$ gives us integrals of the 
form
$$\int_{y \in (\RR^+)^m \times \RR^{l-m}} [dy] P(y) e^{-\langle y,y 
\rangle / 2 \e}$$
where $P(y)$ is a polynomial function of $y \in \RR^l$ and 
$\langle \, ,\, \rangle$ is an inner product on $\RR^l$. Using 
induction on $l$ and calculations similar to those in Example 
\ref{ex9.2}, it follows that $\cali^\epsilon_0(\eta e^{i\baromega}) $ 
is always a polynomial in $\sqrt{\e}$, although not necessarily 
a polynomial in $\epsilon$. 
\end{rem}

\begin{example}
Consider the case we have looked at a number of times 
(Examples 2.4, 3.2, 4.2) 
when $G=SL(2)$ acts on $M = ({\bf P}_1)^{n} $ with $n$ even
(cf. also Section 9 of \cite{JK1}). 
The fixed points of the action of the maximal torus $T$ of $K$ are 
the $n$-tuples of points $(x_1, \dots, x_{n})$ in ${\bf P}_1$ such that 
each $x_j$ is either $0$ or $\infty$, which we index by 
sequences $\underline{\delta} = (\delta_1, \dots, \delta_{n})$ where 
$\delta_j = 1 $ if $x_j = 0 $ and $\delta_j = -1$ if $x_j = \infty$. The 
value of the moment map for such a point is $\sum_{j = 1}^{n} 
\delta_j$. Recall that $H^*(M)$ has $n$ generators of degree $2$, 
and the equivariant cohomology  $H^*_T(M)$ is 
generated by lifts $\xi_i$ ($i = 1, \dots, n$)  of these together with 
one additional generator $\zeta$ of degree 2, subject to the 
relations $(\xi_j)^2 = \zeta^2,$ while $H^*_K(M)$ is generated by 
the $\xi_j$  and $\zeta^2$ subject to the same relations.
Then if $\eta = q (\xi_1, \dots, \xi_{n}, \zeta^2) \in \hk(M)$ we have 
from Theorem 9.1 that
$$\cali^\epsilon_0 (\eta e^{i\baromega})  
= \frac{A_K}{\epsilon^{s/2}} 
( \sum_{
\delta_1 +...+\delta_n >0}
(\int_{-\infty}^{0}dy  e^{- y^2/{2 \epsilon}}y\calp_{-\RR^+} 
 +\int_{0}^{\infty}dy  e^{- y^2/{2 \epsilon}}y\calp_{\RR^+}) $$
$$
 +\sum_{
\delta_1 + ...+ \delta_n =0} \int_{-\infty}^{0}dy  e^{- y^2/{2 \epsilon}}y\calp_{-\RR^+}  
) \times$$
$$  \res_{X=0}^+   \Biggl ( \frac{1}{X^{n-1} } e^{i X(\sum_j \delta_j)}
e^{-i yX } (\prod_j  \delta_j) q(\delta_1 X, \dots, \delta_{n}X, X^2) 
\Biggr ), $$
where $A_K$ is the constant defined in the statement of   Theorem 9.1,
and for a rational function $R(X)$, $\res_{X=0}^+ (e^{i \mu X} R(X) ) $ is the 
coefficient of $1/X$ in the Taylor
expansion if $\mu \ge 0 $ and $0$ otherwise.

An easy calculation shows that the factors $\calp_{\pm \RR^+} \res_{X=0}^+$
 may simply be replaced by $\res_{X=0}$. The contribution of any fixed
point indexed by
$ \{\underline{\delta} \} $ such that $\sum_j \delta_j = 0 $ is 
$$ A_K \epsilon^{-3/2} ~\res_{X=0} \left (  \int_{-\infty}^0 \frac{ (-1)^{n/2} y}{X^{n-1}}
q(\delta_1 X, \dots, \delta_{n} X, X^2) e^{-y^2/{2\epsilon} } e^{- iyX} 
dy \right ). $$
Thus we see from the calculations in Example \ref{ex9.2} that
for some choices of $q$ the individual contributions of
the $F$ with $\mu_T(F)=0$ will involve odd powers of
 $\sqrt{\epsilon}$, although an argument from symmetry shows that
$\cali^\epsilon_0 (\eta e^{i\baromega}) $ is in fact a polynomial
in $\e$ in this example.
\end{example}
\begin{example} As our final example consider the linear
action of $G=\CC^*$ on $M=\PP_n$ with distinct weights 
$r_0,...,r_n \in \ZZ$.
Then $H^*_T(M)$ is generated by two equivariant cohomology
classes $\xi$ and $\zeta$ of degree two subject to the relation
$$\prod_{j=0}^n (\xi - r_j \zeta) =0.$$
The moment map for the action of $T=S^1$ on $M$ is given up
 to a constant
by
$$\mu[x_0,...,x_n]= \frac{r_0|x_0|^2 + ... + r_n|x_n|^2}{|x_0|^2 +
 ... + |x_n|^2}.$$
Suppose now that $r_0=0$ but that $r_j \neq 0$ for $j>0$.Then 
if $\eta = q(\xi,\zeta) \in H^*_T(M)$ we have from Theorem 9.1 that
$$ \cali^\epsilon_0(\eta e^{i\baromega}) = 
\frac{A_K}{\sqrt{\epsilon}}~  \int_{-\infty}^{0}dy 
\res_{X=0}  ( \frac{ e^{-y^2/2\e}e^{-iyX}}{\prod_{k\neq 0}
(r_k X - r_0 X)}q(0,X)  $$
$$+\sum_{j:r_j >0} \int_{- \infty}^{\infty}
dy \res_{X=0}  \frac{e^{ir_j X} e^{-y^2/2\e}e^{-iyX}}{\prod_{k\neq j}
(r_k X - r_j X)}q(r_j X,X) ).$$
The contribution of each $j >0$ to this 
expression is a polynomial in $\e$, but if $q(0,X)=X^N$ then
$$\res_{X=0} \int_{-\infty}^{0}dy \frac{ e^{-y^2/2\e}e^{-iyX}}
{\prod_{k\neq 0}(r_k X - r_0 X)}q(0,X) = \int_{-\infty}^{0}dy 
\frac{ e^{-y^2/2\e}(-iy)^{n-N-2}}{\prod_{k\neq 0}(r_k - r_0 )}$$
which is a nonzero constant multiple of $\e^{(n-N-2)/2}$. Thus if
$n-N$ is odd then $ \cali^\epsilon_0(\eta e^{i\baromega})$ is not
 a polynomial in $\e$.
 \end{example}


\begin{thebibliography}{99}



\bibitem{Aconv} M.F. Atiyah, Convexity and commuting Hamiltonians, 
{\em Bull.London Math. Soc.} {\bf 14} (1982) 1-15.

\bibitem{abym} M.F. Atiyah, R. Bott, The Yang-Mills equations over 
Riemann surfaces. Philos. Trans. Roy. Soc. London Ser. {\bf A 308} 
(1983), no. 1505,523--615.

\bibitem{abmm} M.F. Atiyah, R. Bott, The moment map and 
equivariant cohomology, {\em Topology} {\bf 23} (1984) 1-28.





\bibitem{BBDG} A.A. Beilinson, J. Bernstein, P. Deligne, Faisceaux pervers,
 {\em Ast\'{e}risque} {\bf 100} (1982), 5-171.

\bibitem{BGV} N. Berline, E. Getzler, M. Vergne, {\em Heat Kernels 
 and Dirac Operators}, Springer-Verlag (Grundlehren vol. 298), 1992.

\bibitem{BV1}  N. Berline, M. Vergne, Z\'eros d'un champ de vecteurs
 et classes caract\'eristiques \'equivariantes,
{\em Duke Math. J.} {\bf 50} (1983) 539-549.

\bibitem{BV2} N. Berline, M. Vergne, The equivariant index and 
Kirillov's character formula, {\em Amer. J. Math.} {\bf 107} (1985) 
1159-1190.



\bibitem{BD} T. Br\"ocker, T. Tom Dieck, {\em Representations of 
Compact Lie Groups}, Springer-Verlag, 1985.

\bibitem{cartan} H. Cartan, {Notions d'alg\`ebre diff\'erentielle; 
applications aux vari\'et\'es o\`u op\`ere un groupe de Lie},  
in {\em Colloque de Topologie} (C.B.R.M., Bruxelles, 1950) 15-27; 
La transgression dans un groupe de 
Lie et dans un fibr\'e principal, {\em op. cit.}, 57-71.



\bibitem{DH} J.J Duistermaat, G. Heckman, On the variation in the
cohomology of the symplectic form of the reduced 
phase space, {\em Invent. Math.} {\bf 69} (1982) 259-268; 
Addendum, {\bf 72} (1983) 153-158.



\bibitem{GM1} M. Goresky and R. MacPherson, Intersection homology theory,
{\em Topology } {\bf 19} (1980), 135-162.

\bibitem{GM2} M. Goresky and R. MacPherson, Intersection homology II,
{\em Invent. Math.} {\bf  71} (1983), 77-129.

\bibitem{GH} P. Griffiths, J. Harris, {\em Principles of  algebraic 
geometry,} Wiley-Interscience, 1978.



\bibitem{GK} V. Guillemin, J. Kalkman, The Jeffrey-Kirwan localization 
theorem and residue operations in equivariant cohomology. 
{\em J. Reine Angew.Math.} {\bf 470} (1996), 123--142. 



\bibitem{GSconv} V. Guillemin, S. Sternberg, Convexity properties
of the moment mapping I and II, {\em Invent. Math.}
{\bf 67} (1982) 515-538 and {\bf 77} (1984) 533-546.

\bibitem{GLS} V. Guillemin, E. Lerman, S. Sternberg, On the Kostant 
multiplicity formula, {\em J. Geom. Phys.}
{\bf 5} (1988) 721-750 (1988).

\bibitem{GS} V. Guillemin, S. Sternberg, {\em Symplectic Techniques
in Physics}, Cambridge University Press, 1984.



\bibitem{Hir} H. Hironaka, Resolution of singularities of an algebraic 
variety over a field of characteristic zero. I, II. {\em Ann. of Math.} (2) 
{\bf 79} (1964)109--203; ibid. (2) {\bf 79} 1964 205--326.

\bibitem{JKKW} L.C. Jeffrey, Y.H. Kiem, F.C. Kirwan, J. Woolf,
Intersection pairings on singular moduli spaces of bundles over
a Riemann surface, in preparation.



\bibitem{JK1} L.C.  Jeffrey, F.C. Kirwan, Localization for nonabelian 
group actions, {\em Topology} {\bf 34} (1995) 291-327.



\bibitem{locquant} L.C.  Jeffrey, F.C. Kirwan, Localization and 
the quantization conjecture, {\em Topology} {\bf 36} (1997) 
647-693.



\bibitem{arbrank} L.C.  Jeffrey, F.C. Kirwan, Intersection theory 
on moduli spaces of holomorphic bundles of arbitrary rank on a 
Riemann surface, {\em Ann. Math.} {\bf 148} (1998) 109-196.



\bibitem{Kiem} Y.-H. Kiem, Intersection cohomology of quotients 
of nonsingular varieties, Yale preprint (1999), math.AG/0101254.

\bibitem{KW} Y.-H. Kiem, J. Woolf, The cosupport axiom, equivariant
cohomology and the intersection cohomology of certain symplectic
quotients, preprint math.AG/0101255.

\bibitem{Ki1} F.C. Kirwan, {\em Cohomology of Quotients in Symplectic 
and Algebraic Geometry}, Princeton University Press, 1984.



\bibitem{Ki2} F.C. Kirwan, Partial desingularisations of quotients of
 nonsingular varieties and their Betti numbers, {\em Ann. Math.} 
{\bf 122} (1985), 41-85.

\bibitem{Ki3} F.C. Kirwan, Rational intersection cohomology of 
quotient varieties, {\em Invent. Math} {\bf 86} (1986), 471-505.



\bibitem{luna} D. Luna, Sur les orbites ferm\'{e}s des groupes
alg\'{e}briques r\'eductifs, {\em Invent. Math.} {\bf 16} (1972), 1-5 .

\bibitem{Martin} S.K. Martin, Symplectic                                     
  geometry and gauge theory, Oxford D.Phil. thesis, 1997.

\bibitem{Martin2}S.K. Martin, Symplectic quotients by a 
nonabelian group and by its maximal torus, math.SG/0001002; {\em Annals 
of Mathematics}, to appear.

\bibitem{Martin3} S.K. Martin, Transversality theory, cobordisms, 
and invariants of symplectic quotients, math.SG/0001001;
 {\em Annals of Mathematics},
 to appear.



\bibitem{MSj} E. Meinrenken and R. Sjamaar, Singular reduction
and quantization, {\em Topology} {\bf 38} (1999), 699-762.

\bibitem{Mo} G. Mostow, On a conjecture of Montgomery, {\em
Ann. Math.} {\bf 65} (1957), 513-516.


\bibitem{GIT} D. Mumford, J. Fogarty, F. Kirwan, {\em Geometric 
invariant theory}, 3rd edition, Springer-Verlag, 1994.

\bibitem{New} P.E. Newstead,  {\em Introduction to moduli problems
and orbit spaces}, Tata Institute Lectures 51, Springer-Verlag, 1978.

\bibitem{paradan} P.-E. Paradan, The moment map and equivariant
cohomology with generalized coefficients, {\em Topology} {\bf 39} (2000),
401-444.





\bibitem{TW} S. Tolman and J. Weitsman,  The cohomology rings
of abelian symplectic quotients, preprint math.DG/9807173.





\bibitem{V2} M. Vergne, A note on the Jeffrey-Kirwan-Witten 
localisation formula. {\em Topology} {\bf 35}  (1996),  243--266. 


\bibitem{tdgr} E.  Witten, { Two dimensional gauge theories
 revisited}, preprint hep-th/9204083; {\em J. Geom. Phys.} {\bf 9} 
(1992) 303-368.



\bibitem{Woolf} J. Woolf, {\em Some Topological Invariants
of Singular Symplectic Quotients}, D.Phil. thesis, Oxford, 1999.

\bibitem{Woolf2} J. Woolf, The decomposition theorem and the 
intersection cohomology of quotients in algebraic geometry, preprint 
math.AG/0110137
(2000).





\end{thebibliography}
\end{document}